\documentclass[a4paper,10pt]{article}

\usepackage{amsmath,amssymb,amsthm}

\def\ukn{\mathbf{U}_{\mathbb{K},n}^-}

\def\uan{\mathbf{U}_{n}^-}

\def\uasn{\mathbf{U}^-(\hat{\mathfrak{sl}}_n)}
\def\uksn{\mathbf{U}_\K^{-}(\hat{\mathfrak{sl}}_n)}

\def\La{\Lambda^\infty}

\def\C{\mathbb{C}}
\def\Z{\mathbb{Z}}
\def\N{\mathbb{N}}
\def\A{\mathbb{A}}
\def\K{\mathbb{K}}
\def\S{{\mathbb{S}}}
\def\R{{\mathbf{R}}}
\def\vac{|0\rangle}
\def\qed{$\hfill \blacksquare$}

\def\td{\tilde{d}}
\def\m{\mathbf{m}}
\def\f{\mathbf{f}}
\def\Si{\mathfrak{S}}
\def\ASi{\widehat{\mathfrak{S}}}
\def\AH{\widehat{\mathbf{H}}}
\def\x{\text{x}}
\def\i{\mathbf{i}}
\def\TD{\bigotimes^D}
\def\U{\mathbf{U}}

\def\s{\tau}
\def\b{\mathbf{b}}
\def\d{\mathbf{d}}
\def\e{\mathbf{e}}
\def\k{\mathbf{k}}
\def\s{\mathbf{s}}

\newtheorem*{theo}{Theorem}
\newtheorem*{prop}{Proposition}
\newtheorem{lem}{Lemma}
\newtheorem*{cor}{Corollary}
\numberwithin{equation}{section}
\numberwithin{lem}{section}

\begin{document}
\title{The Hall algebra of a cyclic quiver\\
and canonical bases of Fock spaces}
\author{Olivier Schiffmann}
\date{}
\maketitle
\noindent

\paragraph{0.1 Introduction.} The Lusztig formula for characters of the
simple finite-dimensional modules of the quantum group
$\U_\epsilon(\mathfrak{sl}_k)$, where
$\epsilon=e^{2i\pi/n}$, gives the multiplicity of
the Weyl module $W_\mu$ of $\U_\epsilon(\mathfrak{sl}_k)$ with highest weight $\mu$
in the simple $\U_\epsilon(\mathfrak{sl}_k)$-module $V_\lambda$ with highest weight
$\lambda$. Namely, 
$$[V_\lambda:W_\mu]=\sum_{y}(-1)^{l(xy)} P_{y,x}(1),$$
where $x \in \widehat{\Si}_k$ is minimal such that $\nu=\lambda.x^{-1}$
satisfies $\nu_i<\nu_{i+1}$ for $i=1,2 \ldots k-1$ and $\nu_i-\nu_k
\geq 1-k-n$, and $\mu=\lambda. x^{-1}y$. This conjecture is proved by
Kazhdan-Lusztig \cite{KL} and Kashiwara-Tanisaki \cite{KT}. The proof
relies on an equivalence between the category of
finite-dimensional $\U_\epsilon(\mathfrak{sl}_k)$-modules and a category of negative-level
representations of the affine algebra $\hat{\mathfrak{sl}}_k$ which are integrable with respect to $\mathfrak{sl}_k$. In
\cite{VV}, Varagnolo and Vasserot propose a new approach to this
conjecture, based on the geometric constructions of simple
finite-dimensional $\U_q(\hat{\mathfrak{sl}}_k)$-modules of \cite{GV}
and on the theory of canonical bases. Let $\U^-_n$ be the generic Hall
algebra of the cyclic quiver of type $A_{n-1}^{(1)}$ (defined over the
ring $\C[q,q^{-1}]$) and let $\mathbf{B}$ be the intersection
cohomology basis of $\U^-_n$. Let $\Lambda^\infty$ be the Fock space
representation of $\U^-_n$ (see \cite{KMS} and \cite{VV}), and
let $\mathbf{B}^{\pm}$ be the Leclerc-Thibon canonical bases of
$\Lambda^\infty$ (see \cite{LT}). Varagnolo and Vasserot show that the
Lusztig formula follows from the equality
$\mathbf{B}\vac=\mathbf{B}^+$, where $\vac$ is the vacuum
vector of $\Lambda^\infty$. In this paper, we give a direct proof of
this equality, which can be thought of as
a q-analogue of the Lusztig conjecture. This also yields a proof of
the positivity conjecture for the basis $\mathbf{B}^+$ (see \cite{LLT1}, Conjecture 6.9 i)). Note that the equality
$\mathbf{B}\vac=\mathbf{B}^+$ does not follow from the general theory
developped in \cite{K} or \cite{Lb} since the Fock space is not
irreducible as a $\U_q^-(\hat{\mathfrak{sl}}_n)$-module. We also formulate and prove an analogue of the conjecture of Varagnolo and Vasserot for higher level Fock space representations $\Lambda^{\infty}_{\mathbf{s}_l}$ of $\U^-_n$ (see \cite{JMMO}). The corresponding canonical bases have been recently introduced by Uglov, \cite{U}.\\
\hbox to1em{\hfill}Our proof relies on the construction of an isomorphism
$$\U^-_n \simeq \U_q^-(\hat{\mathfrak{sl}}_n) \otimes
\C[\mathbf{z}_1,\mathbf{z}_2,\ldots].$$
The central subalgebra $\C[v,v^{-1}][\mathbf{z}_1,\mathbf{z}_2,\ldots]$ of $\U^-_n$ is linearly spanned by elements of the \textit{dual} canonical basis of $\U^-_n$. By a quantum affine version of the Schur-Weyl duality, it also admits an interpretation in terms of the Bernstein center $Z(\AH_l)$ of the affine Hecke algebra $\AH_l$ of type $A_l$ for $l \leq n$. We consider a basis $(a_\lambda)$ of the
central subalgebra
$\C[v,v^{-1}][\mathbf{z}_1,\mathbf{z}_2,\ldots]$ of $\U^-_n$ corresponding to Schur polynomials
and describe its action on $\Lambda^\infty$. The equality
$\mathbf{B}\vac=\mathbf{B}^+$ follows from the caracterization of
$\mathbf{B}^+$ in terms of a lattice and the $\bar{\ }$-involution of 
$\Lambda^\infty$.\\
\hbox to1em{\hfill}Finally, we note that a proof of the Lusztig conjecture (but not of its q-analogue) using the above ideas has been recently found by B. Leclerc (see \cite{Lec}).

\paragraph{0.2 Notations.} Set $\S=\C[v],\;\A=\C[v,v^{-1}]$ and
$\K=\C(v)$. Throughout the paper we fix some integer $n>1$. Let $q$ be
a prime power and let $\mathbb{F}$ be a finite field with $q^2$
elements. Let $(\epsilon_i)$, $i \in \Z/n\Z$ be the canonical basis of
$\N^{\Z/n\Z}$. For $i \in \Z/n\Z$ and $l \in \N^*$, define the \textit{cyclic segment} $[i;l)$ to be the image of the projection to $\Z/n\Z$ of the segment $[i_0,i_0+l-1] \subset \Z$ for any $i_0 \equiv i\;(\mathrm{mod}\;n)$. A {\textit{cyclic multisegment} is a linear combination $\m=\sum_{i,l} a^l_i [i;l)$ of cyclic segments with coefficients $a^l_i \in \N$. Let $\mathcal{M}$ be the set of cyclic multisegments. For $\m \in \mathcal{M}$, we put $|\m|= \sum_{i,l}
a^l_i$ and $\mathrm{dim}\;\m=\sum_{l,i}a^l_i(\epsilon_{i}+\ldots + \epsilon_{i+l-1})  \in \N^{\Z/n\Z}$. Let $\Pi$ denote the set of partitions and let $\Pi_D \subset \Pi$ be the subset formed by partitions with at most $D$ parts.
\noindent
\section{Hall algebra of the cyclic quiver.}
\paragraph{1.1} Let $Q$ be the quiver of type $A^{(1)}_{n-1}$, i.e the oriented graph with vertex set $I=\Z/n\Z$ and edge set $\Omega=\{({i},{i+1}), i \in I\}$. For any $I$-graded $\mathbb{F}$-vector space $V=\bigoplus_{i \in I} V_i$ let $E_V \subset \bigoplus_{(i,j) \in \Omega} \mathrm{Hom}\;(V_i,V_j)$ denote the space of nilpotent representations of $Q$, i.e collections of linear maps $(\sigma_i: V_i \to V_j, (i,j) \in \Omega)$ satisfying the following condition: for any $i_1 \in I$ there exists $r \in \N$ such that $\sigma_{i_r}\ldots \sigma_{i_2}\sigma_{i_1}=0$ for any $i_2,\ldots i_r$. The group $G_V=\prod_{i\in I} GL(V_i)$ acts on $E_V$ by conjugation. For each $i \in I$ there exists a unique simple $Q$-module $S_i$ of dimension $\epsilon_i$, and for each pair $(i,l) \in I \times \N^*$ there exists a unique (up to isomorphism) indecomposable $Q$-module $S_{i;l}$ of length $l$ and tail $S_i$. Furthermore, every nilpotent $Q$-module $M$ admits an essentially unique decomposition
\begin{equation}\label{E:1.1}
M \simeq \bigoplus_{i,l} a_i^lS_{i;l}.
\end{equation}
The classification of $Q$-modules is independent of the base field. We denote by $\overline{\mathbf{m}}$ the isomorphism class of $Q$-modules corresponding (by (\ref{E:1.1})) to the multisegment $\mathbf{m}=\sum_{i,l} a_i^l [i;l)$. 
\paragraph{1.2} We recall the Lusztig construction (see \cite{Lb}). For any $I$-graded vector space $V$, let $\C_{G}(E_V)$ be the set of $G_V$-invariant functions $E_V \to \C$. For each $\mathbf{d} \in \N^I$, let us fix an $I$-graded vector space $V_{\mathbf{d}}$ of dimension $\mathbf{d}$. Given $\mathbf{a},\mathbf{b} \in \N^I$ such that $\mathbf{a}+\mathbf{b}=\mathbf{d}$, consider the diagram
$$E_{V_\mathbf{a}}\times E_{V_\mathbf{b}} \stackrel{p_1}{\leftarrow} E \stackrel{p_2}{\to}F\stackrel{p_3}{\to} E_{V_{\mathbf{d}}},$$
where $E$ is the set of triples $(x,\phi,\psi)$ such that $x \in E_{V_{\mathbf{d}}}$, 
$$0 \to V_{\mathbf{a}} \stackrel{\phi}{\to} V_{\mathbf{d}} \stackrel{\psi}{\to} V_{\mathbf{b}}\to 0$$
is an exact sequence of $I$-graded vector spaces and $\phi(V_{\mathbf{a}})$ is stable by $x$, and where $F$ is the set of pairs $(x,U)$ such that $U \subset V_{\mathbf{d}}$ is an $I$-graded, $x$-stable subspace of dimension $\mathbf{a}$.\\
\hbox to1em{\hfill}Given $f \in \C_{G}(E_{V_{\mathbf{a}}})$ and $g \in \C_G(E_{V_{\mathbf{b}}})$, set
\begin{equation}\label{E:1.2} 
f \circ g= q^{-m(\mathbf{b},\mathbf{a})} (p_3)_{!}h \in \C_G(V_{\mathbf{d}}),
\end{equation}
where $h\in \C(F)$ is the unique function satisfying $p_2^*(h)=p_1^*(fg)$ and $m(\mathbf{b},\mathbf{a})=\sum_{i\in I} a_ib_i + \sum_{(i_1,i_2) \in \Omega} b_{i_1}a_{i_2}$. Set $\U^-_{q,n} = \bigoplus_{\mathbf{d}} \C_G(E_{V_\mathbf{d}})$. Then $(\U^-_{q,n}, \circ)$ is an associative algebra. \\
\hbox to1em{\hfill}For $\m \in \mathcal{M}$ with
$\mathrm{dim}\;\m=\mathbf{d}$, we let $O_\m \subset E_{V_\mathbf{d}}$ be the $G_{V_\mathbf{d}}$-orbit consisting of representations in the class $\overline{\mathbf{m}}$, and we let $\mathbf{1}_\mathbf{m} \in \C_G(V_{\mathbf{d}})$ be the characteristic function of $O_\m$. Finally, we set $\f_\m=q^{-\mathrm{dim}\;O_\m} \mathbf{1}_\m$, and if $\mathbf{d} \in \N^I$, we let $\f_\mathbf{d}$ be the characteristic function of the trivial representation in $V_\mathbf{d}$. Note that, by definition, $(\f_\m)_{\m \in \mathcal{M}}$ is a $\C$-basis of $\U^-_{q,n}$.
\paragraph{1.3} Let $\mathbf{a},\mathbf{b},\mathbf{d} \in \N^I$ such that $\mathbf{a}+\mathbf{b}=\mathbf{d}$. Fix a subspace $U$ of $V_{\mathbf{d}}$ of dimension $\mathbf{a}$ and a pair of graded vector apace isomorphisms $U \simeq V_\mathbf{a},\;V_\mathbf{d}/U \simeq V_\mathbf{b}$. Consider the diagram
$$E_U \times E_{V/U} \stackrel{p}{\leftarrow} E \stackrel{i}{\to} E_{V_\mathbf{d}},$$
where $E \subset E_{V_\mathbf{d}}$ is the subset of the representations preserving $U$. Set
$$
\Delta_{\mathbf{a},\mathbf{b}}:\;\C_G(V_\mathbf{d})  \to \C_G(V_\mathbf{a}) \otimes \C_G(\mathbf{b}),\qquad
f  \mapsto q^{-n(\mathbf{b},\mathbf{a})} p_{!}i^*f.
$$
Here $n(\mathbf{b},\mathbf{a})=\sum_{(i_1,i_2) \in \Omega} b_{i_1}a_{i_2} -\sum_{i \in I}a_ib_i$.
\paragraph{1.4} It is known
that the structure constants of $(\U_{q,n}^-,\circ,
(\Delta_{\mathbf{a},\mathbf{b}}))$ in the basis $(\f_\m)$ are values
at $q=v^{-1}$ of some universal Laurent polynomials in $\A$. Thus there exists an algebra $\uan$ defined over $\A$ of which $\U^-_{q,n}$ is the specialisation at $v=q^{-1}$. This algebra is called the \textit{generic Hall algebra}. The algebra $\uan$ is naturally
$\N^I$-graded and we will denote by $\uan[\alpha]$ the graded
component of degree $\alpha \in \N^I$. Let $\U_\K(\hat{\mathfrak{sl}}_n)$ (resp. $\U(\hat{\mathfrak{sl}}_n)$) denote the rational form (resp. the Lusztig integral form) of the quantum affine algebra of type $A_{n-1}^{(1)}$, and let $e_i, k_i,f_i$, $i\in I$ (resp. $e_i^{(l)},k_i,f_i^{(l)}$, $i\in I$, $l \in \N$) be the standard Chevalley generators (resp. their quantized divided powers). Let $\uasn$ be the subalgebra of $\U(\hat{\mathfrak{sl}}_n)$ generated by $f^{(l)}_i, i\in I$, $l \in \N^*$. It is known that the map $f^{(l)}_i \mapsto
\mathbf{f}_{l\epsilon_i}$ extends to an embedding of the algebras $\uasn \hookrightarrow
\uan$. Set $\ukn=\uan \otimes_\A \K$.

\paragraph{1.5} Let $\U^0$ be the commutative
$\A$-algebra generated by elements $\mathbf{k}_\mathbf{d},\;\mathbf{d} \in \Z^I$
satisfying
$$\mathbf{k}_0=1,\qquad
\mathbf{k}_{\mathbf{a}}\mathbf{k_b}=\mathbf{k}_{\mathbf{a}+\mathbf{b}},
\qquad \forall\; \mathbf{a},\mathbf{b} \in \Z^I.$$
Set $\widetilde{\U}^-_n=\uan \otimes_\A \U^0_{}$ and put
$$(f \otimes \mathbf{k}_{\mathbf{a}}) \circ (g \otimes
\mathbf{k}_\mathbf{b})=v^{\mathbf{a\cdot d}} (f \circ g)\otimes
(\mathbf{k}_{\mathbf{a} +\mathbf{b}}), \qquad \forall\;g \in
\uan[\mathbf{d}],\; \forall\;f \in \uan,$$
where $\mathbf{a \cdot
  d}=n(\mathbf{a},\mathbf{d})+n(\mathbf{d},\mathbf{a})$. Finally, define
$$
\tilde{\Delta}: \widetilde{\U}_n^- \to  \widetilde{\U}_n^- \otimes
 \widetilde{\U}_n^- ,\qquad
f \otimes \mathbf{k}_\mathbf{c} \mapsto \sum_{\mathbf{d}=\mathbf{a} +
 \mathbf{b}}
 \Delta_{\mathbf{a},\mathbf{b}}(f)(\mathbf{k}_{\mathbf{b}+\mathbf{c}}
 \otimes \mathbf{k}_\mathbf{c})\qquad \forall f \in \uan[\mathbf{d}].
$$
It is proved in \cite{G} that $(\widetilde{\U}_n^-,\circ,
\tilde{\Delta})$ is a bialgebra.

\paragraph{1.6} Define the following symmetric bilinear form on
$\widetilde{\U}^-_n$ :
$$\langle \f_\m \otimes \mathbf{k}_{\mathbf{a}},\f_{\m'}\otimes
  \mathbf{k}_{\mathbf{b}}\rangle=v^{-(\mathbf{a}+\mathbf{b})\mathbf{\cdot d}-2\; \mathrm{dim}\;\mathrm{Aut}(\m)} \frac{(1-v^2)^{|\m|}}{|\mathrm{Aut}(\m)|} \delta_{\m,\m'},$$
where $\mathbf{d}=\mathrm{dim}\;\m$ and $\mathrm{Aut}(\m)$ stands for
  the group of automorphism of any representation in the orbit
  $O_\m$. For any $f,g,h \in \widetilde{\U}^-_n$ we have (see \cite{G})
$$\langle fg,h \rangle=\langle f \otimes g,\tilde{\Delta}(h)\rangle.$$
It is clear that the restriction of $\langle\,,\,\rangle$ to the subalgebra $\uan$
is nondegenerate.

\paragraph{1.7} For $i \in I$, let $e'_i: \uan \to \uan$ be the adjoint
of the left multiplication by $\f_i$. It is a homogeneous operator of
degree $-\epsilon_i$. Let $\tilde{f}_i,\tilde{e}_i$ be the
Kashiwara operators (see \cite{K}, Section 3). Recall that a crystal
basis of $\uan$ is a pair $(\mathcal{L},\text{B})$ where $B\subset
\mathcal{L}/v\mathcal{L}$ is a $\C$-basis satisfying the following
conditions :
\begin{enumerate}
\item[i)] for any $i \in I$ we have $\tilde{e}_i\mathcal{L},\tilde{f}_i \mathcal{L}
  \subset \mathcal{L}$, and $\tilde{e}_i(\text{B}),\tilde{f}_i(\text{B}) \subset \text{B} \cup
  \{0\}$,
\item[ii)] for any $i \in I, \text{b},\text{b}'\in \text{B}$ we have $\tilde{e}_i(\text{b})=\text{b}'$
  if and only if $\tilde{f}_i(\text{b}')=\text{b}$.
\end{enumerate}
 Set $\mathcal{L}=\bigoplus_\m \S \;\f_\m$. Let
$\text{b}_\m$ denote the class of $\f_\m$ in
$\mathcal{L}/v\mathcal{L}$ and set $\text{B}=\{\text{b}_\m,\;\m \in
\mathcal{M}\}$. The following result is proved in \cite{LTV},
Theorem. 4.1 :
\begin{theo} The couple $(\mathcal{L},B)$ is a crystal basis of
  $\uan$.
\end{theo}
\noindent The crystal graph $\mathcal{C}$ of $\uan$ has the vertex set $\mathcal{M}$ and edges
$\m \stackrel{i}{\to} \m'$ whenever $\tilde{f}_i(\text{b}_\m)=\text{b}_{\m'}$. It is
explicitely described in \cite{LTV}. Call a multisegment $\m=\sum_{i,l} a^l_i[i;l)$ \textit{completely periodic} if $a^l_i=a^l_j$ for all $l \in \N^*$ and $i,j \in I$. Let $\mathcal{M}^{\mathrm{per}}$ be the
set of all completely periodic multisegments. Then $\tilde{e}_i(\text{b}_\m)=0$ for each $\m \in \mathcal{M}^{\mathrm{per}}$ and $i \in I$, and the connected component of $\mathcal{C}$ containing $\text{b}_\m$ is isomorphic to the crystal graph of $\uasn$.

\paragraph{1.8} For $\mathbf{m}\in \mathcal{M}$, set
$$\mathbf{b}_{\mathbf{m}}=\sum_{i,\mathbf{n}}v^{-i+\mathrm{dim}\,O_{\mathbf{m}}+\mathrm{dim}\,O_{\mathbf{n}}}\mathrm{dim}\,\mathcal{H}^i_{O_{\mathbf{n}}}(IC_{O_{\mathbf{m}}}) \mathbf{f}_{\mathbf{n}},$$
where $\mathcal{H}^i_{O_{\mathbf{n}}}(IC_{O_{\mathbf{m}}})$ is the
stalk over a point of $O_{\mathbf{n}}$ of the ith intersection
cohomology sheaf of the closure $\overline{O}_{\mathbf{m}}$ of
$O_{\mathbf{m}}$. Then $\mathbf{B}=\{\mathbf{b}_{\mathbf{m}}\}$ is the
canonical basis of $\uan$, introduced in \cite{VV}. There exists a unique semilinear ring involution $x
\mapsto \overline{x}$ of
$\uan$ satisfying $\overline{\mathbf{b}_\m}=\mathbf{b}_\m$ (see \cite{VV}, Proposition 7.5). The element $\mathbf{b}_{\mathbf{m}}$ is caracterized
by the following two properties :
$$ \text{i)}\;\overline{\mathbf{b}_{\mathbf{m}}}=\mathbf{b}_{\mathbf{m}},
\qquad \mathrm{and}\qquad  \text{ii)}\;
\mathbf{b}_{\mathbf{m}} \in \mathbf{f}_{\mathbf{m}}+  v
\mathcal{L}.$$
 
\hbox to1em{\hfill}Call a multisegment $\m=\sum_{i,l} a^l_i$ \textit{aperiodic} if
for each $l \in \N^*$ there exists $i \in I$ such that $a^l_i=0$, and let $\mathcal{M}^{\mathrm{ap}}$ be the set of aperiodic
multisegments. Lusztig proved in \cite{L2} that $\{\mathbf{b}_\m\;|\m \in \mathcal{M}^{\mathrm{ap}}\}$ is the global canonical basis of $\uasn$.
\paragraph{1.9} Fix $D \in \N$. Let $\Si_D$ (resp. $\ASi_D$, resp. $\AH_D$) be the symmetric group (resp. affine symmetric group, resp. affine Hecke algebra) of type $GL_D$. The $\A$-algebra $\AH_D$ is generated by elements $T_i^{\pm}, X_j^{\pm}, \; i \in [1,D-1], \;j \in [1,D]$ with relations
\begin{alignat*}{2}
&T_i\,T_i^{-1}=1=T_i^{-1}\,T_i,\qquad & \qquad &(T_i+1)(T_i-v^{-2})=0,\\
&T_i\,T_{i+1}\,T_i=T_{i+1}\,T_i\,T_{i+1},\qquad&\qquad&
|i-j|>1\Rightarrow T_i\,T_j=T_j\,T_i,\\
&X_i\,X_i^{-1}=1=X_i^{-1}\,X_i,\qquad&\qquad &X_i\,X_j=X_j\,X_i,\qquad\\
&T_i\,X_i\,T_i=v^{-2}X_{i+1},\qquad&\qquad 
&j\not= i,i+1\Rightarrow X_j\,T_i=T_i\,X_j.
\end{alignat*}
Let $\AH_\infty$ be the $\A$-algebra generated by $T_i^{\pm},X_i^{\pm},\;i,j \in \N^*$ with the same relations as above.\\
\hbox to1em{\hfill}The center of $\AH_D$ is $Z(\AH_D)=\A[X_1^{\pm 1},\ldots ,X_D^{\pm 1}]^{\Si_D}$. Set 
$$Z^+_D=\A[X_1,\ldots ,X_D]^{\Si_D}, \qquad Z^-_D=\A[X_1^{-1},\ldots ,X_D^{-1}]^{\Si_D}.$$

\paragraph{1.10} Let $\A^{(\Z)}$ be the $\A$-linear span of vectors $\x_i$, $i \in \Z$. Following \cite{VV}, Section 8.1, let $\widetilde{\U}^-_n$ act on $\A^{(\Z)}$ by
\begin{align}
\f_\m(\x_i)&=\sum_{j \geq i} \delta_{\m,[\overline{i},j-i+1)} \x_{j+1}, \qquad \forall\; \m \in \mathcal{M},\label{E:1102}\\
\mathbf{k}_{\alpha}(\x_i)&=v^{-n(\alpha,\epsilon_{\overline{i}})}\x_i\qquad \forall \alpha \in \N^I.\label{E:1103}
\end{align}
Set $\bigotimes^D=(\A^{(\Z)})^{\otimes D}$ and let $\widetilde{\U}^-_n$ act on $\bigotimes^D$ via the coproduct $\tilde{\Delta}$. For $\i=(i_1,\ldots i_D) \in \Z^D$ we set $\otimes \x_\i=\x_{i_1} \otimes \ldots \otimes \x_{i_D} \in \TD$. Then $\AH_D$ acts on $\TD$ on the right in the following way
\begin{align}
&(\otimes \x_\i)T_k=\cases v^{-2}\otimes \x_\i & \mathrm{if}\; i_k=i_{k+1}\\
v^{-1}\otimes \x_{(\i)s_k} & \mathrm{if}\; -n<i_k<i_{k+1} \leq 0\\
v^{-1}\otimes \x_{(\i)s_k} + (v^{-2}-1)\otimes \x_\i & \mathrm{if} -n<i_{k+1} \leq 0, \endcases \label{E:110}\\
&(\otimes \x_\i)X_{j}=\otimes \x_{(\i-n\epsilon_j)},\label{E:1101}
\end{align}
where $k \in [1,D-1], \; j \in [1,D]$ and $s_k \in \Si_D$ is the $k$th simple transposition. Moreover, the actions of $\widetilde{\U}^-_n$ and $\AH_D$ on $\TD$ commute (see \cite{VV}, Section 8.2).

\paragraph{1.11} Let $\Omega^D=\sum_i \mathrm{Im}\;(1+T_i) \subset \TD$. For any $\i \in \Z^D$ let $\wedge \x_\i$ be the class of $\otimes \x_\i$ in the quotient $\TD/\Omega^D$. Then
$$\{\wedge \x_\i \;|\;i_1>i_2\ldots > i_D\}$$
is a basis of $\TD/\Omega^D$ (see \cite{KMS}, Proposition 1.3). If $\lambda=(\lambda_1 \geq \ldots \geq \lambda_D) \in \Pi_D$ set $|\lambda \rangle =\wedge \x_\i$ where $i_k=\lambda_k+1-k$, $k=1,\ldots D$. Let $\Lambda^D$ be the $\A$-linear span of the vectors $|\lambda \rangle$, $\lambda \in \Pi_D$. The representation of $\widetilde{\U}^-_n$ on $\TD$ descends to $\TD/\Omega^D$ and restricts to $\Lambda^D$ (see \cite{VV}, 9.2). Note that $Z(\AH_D)$ acts on $\TD/\Omega^D$ and that $Z^-_D$ acts on $\Lambda^D$.\\
\hbox to1em{\hfill}Let $\bigotimes^\infty$ be the $\A$-linear span of semi-infinite monomials
$$\otimes \x_\i=\x_{i_1} \otimes \x_{i_2} \otimes \ldots$$
such that $i_k=1-k$ for $k \gg 0$. The algebra $\AH_\infty$ acts on $\bigotimes^\infty$ via formulas (\ref{E:110}), (\ref{E:1101}). Let $\Omega^\infty =\sum_i \mathrm{Im}\;(1+T_i) \subset \bigotimes^\infty$. As before, if $\i=(i_1,i_2,\ldots)$ is a sequence of integers such that $i_k=1-k$ for $k \gg 0$, let $\wedge \x_\i$ denote the class of $\otimes \x_\i$ in the quotient $\bigotimes^\infty/\Omega^\infty$. If $\lambda \in \Pi$ we set $|\lambda \rangle=\wedge\x_\i$ where $i_k=\lambda_k+1-k$. Finally, set $\Lambda^\infty=\bigoplus_\lambda \A |\lambda \rangle$. It is shown in \cite{VV}, Section 10.1 that $\uan$ acts on $\Lambda^\infty$. Its restriction to $\uasn$ coincides with the level 1 integrable module $F_{(0)}$ considered in \cite{KMS} (see also \cite{H}).
\noindent
\section{A central subalgebra of $\uan$.}
\paragraph{2.1} Set $\R=\bigcap_i \mathrm{Ker}\;e'_i=\bigcap_i
\mathrm{Ker}\;\tilde{e}_i\subset \uan$ and put $\R_\K=\R \otimes_\A \K$. For simplicity, we set $\s=(1,1\ldots 1) \in \N^I$. We first show
\begin{prop} The following properties hold.
\begin{enumerate}
\item[\text{i)}] $\R$ is a graded subalgebra of $\uan$ satisfying $\overline{\R}=\R$ and 
\begin{equation}\label{E:22}
\mathrm{dim}\;\R_\K[\alpha]=\cases 0 & \mathrm{if}\;\alpha \not\in \N \s \\ p(k)&\mathrm{if}\;\alpha=k\s. \endcases
\end{equation}
where $p(k)$ is the number of partitions of the integer $k$.
\item[\text{ii)}] The subalgebras $\R$ and $\uasn$ commute and the multiplication map $m$ induces an isomorphism
$$m: \uasn \otimes_\A \R \stackrel{\sim}{\to} \uan.$$
\item[\text{iii)}] The pair $(\R,\Delta)$ is a bialgebra.
\end{enumerate}
\end{prop}
\noindent
\textit{Proof.} It is clear that $\R$ is graded. Moreover, it follows from \cite{LTV}, Theorem 4.1, that
$$\mathrm{dim}\;\left(\R_\K[\alpha]\cap \mathcal{L}/v(\R_\K[\alpha]\cap \mathcal{L})\right)=\cases 0 & \mathrm{if}\;\alpha \not\in \N \s \\ p(k)&\mathrm{if}\;\alpha=k\s. \endcases
$$
Hence $\R_\K[\alpha]\cap \mathcal{L}$ is a free $\S$-module of the given dimension, and (\ref{E:22}) follows. That $\R$ is a subalgebra is a consequence of the following equality
\begin{equation}\label{E:221}
e'_i(ab)=v^{-\mathrm{wt}_i(\alpha)}ae'_i(b)+e'_i(a)b,\qquad \forall\;a \in \R[\alpha],\;b \in \R[\beta]
\end{equation}
where we set $\mathrm{wt}_i(\gamma)=2\gamma_i-\gamma_{i-1}-\gamma_{i+1}$ for all $\gamma \in \N^I$. To prove (\ref{E:221}), observe that for any $c \in \uan$, we have
\begin{equation*}
\begin{split}
\langle e'_i(ab),c\rangle=\langle ab, \f_ic \rangle&=\langle a \otimes b,(\f_i \otimes 1+ \mathbf{k}_{\epsilon_i}\otimes \f_i)\Delta(c)\rangle\\
&=\sum \langle a, \f_i c'\rangle \langle b,c''\rangle +\langle a,\mathbf{k}_{\epsilon_i}c'\rangle \langle b,\f_i c''\rangle \\
&=\langle e'_i(a)b+v^{-\mathrm{wt}_i(\alpha)}ae'_i(b),c\rangle,
\end{split}
\end{equation*}
where we use Sweedler's notation $\Delta(c)=\sum c'\otimes c''$.\\
\hbox to1em{\hfill}We now prove the last statement in $i)$. For $i \in I$, consider $e''_i: \uan \to \uan,\; x \mapsto \overline{e'_i(\overline{x})}$. We claim that, for all $i, j \in I$,
\begin{equation}\label{E:222}
e'_ie''_j=v^{\mathrm{wt}_i(\epsilon_j)}e''_je'_i.
\end{equation}
To prove (\ref{E:222}), set $S=e'_ie''_j-v^{\mathrm{wt}_i(\epsilon_j)}e''_je'_i$. Then
\begin{align*}
e'_ie''_j\mathbf{f}_{l\s}=&e'_i(\mathbf{f}_{l\s}e''_j+v^{-1} \mathbf{f}_{l\s-\epsilon_j})\\
=& \mathbf{f}_{l\s}e'_ie''_j+v \mathbf{f}_{l\s-\epsilon_i}e''_j + v^{\mathrm{wt}_i(\epsilon_j)-1}\mathbf{f}_{l\s-\epsilon_j}e'_i+v^{\delta_{ij}-\delta_{i-1,j}}\mathbf{f}_{l\s-\epsilon_i-\epsilon_j}.
\end{align*}
Similarly,
\begin{align*}
v^{\mathrm{wt}_i(\epsilon_j)}&e'_ie''_j\mathbf{f}_{l\s}\\
=&v^{\mathrm{wt}_i(\epsilon_j)}\big\{ \mathbf{f}_{l\s}e''_je'_i+v^{-1} \mathbf{f}_{l\s-\epsilon_j}e'_i + v^{1-\mathrm{wt}_i(\epsilon_j)}\mathbf{f}_{l\s-\epsilon_i}e''_j+v^{-\delta_{ij}+\delta_{i+1,j}}\mathbf{f}_{l\s-\epsilon_i-\epsilon_j}\big\}\\
=&\mathbf{f}_{l\s}v^{\mathrm{wt}_i(\epsilon_j)}e''_je'_i+v^{\mathrm{wt}_i(\epsilon_j)-1} \mathbf{f}_{l\s-\epsilon_j}e'_i + v\mathbf{f}_{l\s-\epsilon_i}e''_j+ v^{\delta_{ij}-\delta_{i-1,j}}\mathbf{f}_{l\s-\epsilon_i-\epsilon_j}.
\end{align*}
Hence $S\f_{l\s}=\f_{l\s}S$ for all $l \in \N^*$. In an analogous fashion, $S\f_k=\f_kS$ for all $k \in I$. 
Moreover, by \cite{GP}, the algebra $\ukn$ is generated by $\f_k,\;k \in I$ and $\f_{l\s}\; l \in \N^*$. Hence $S=0$ and (\ref{E:222}) is proved.\\
\hbox to1em{\hfill}We now turn to $ii)$. By (\ref{E:222}) we have $e''_j(\mathrm{Ker}\;e'_i) \subset \mathrm{Ker}\;e'_i$. Hence $e''_j(\R) \subset \R$. But $e''_j$ is a homogeneous operator of degree $-\epsilon_j$, and by (\ref{E:22}) the only non-zero graded components of $\R$ are located in degrees $ls,\;l \in \N^*$. Thus $e''_j(\R)=0$ for all $j$, and $\R \subset \bigcap_j \mathrm{Ker}\;e''_j=\overline{\R}$. Hence $\overline{\R}=\R$ as desired.\\
\hbox to1em{\hfill}By (\ref{E:221}) we have
$$e'_i(\f_jx)=\delta_{ij}x=e'_i(x\f_j),\;\qquad \forall x \in \R,\;j \in I.$$
Hence $\R$ is $\mathrm{ad}\;\f_j$-stable. But $\mathrm{ad}\;\f_j$ is
homogeneous of degree $\epsilon_j$. This again implies that
$(\mathrm{ad}\;\f_j)_{|\R}=0$. Thus the subalgebras $\uasn$ and
$\R$ commute.\\
\hbox to1em{\hfill}The operators $e'_i$ are locally nilpotent. For any $u \in \uan$ there exists a sequence $i_1,i_2, \ldots i_r$ such that $e'_{i_r} \ldots e'_{i_1}u \in \R$. Since the operators $e'_i$ and $\tilde{e}_i$ are proportional, it follows that $\tilde{e}_{i_r}\ldots \tilde{e}_{i_1}u \in \R$. Then $u=\tilde{f}_{i_1}\ldots \tilde{f}_{i_r} \tilde{e}_{i_r}\ldots \tilde{e}_{i_1}u$. Therefore the multiplication map $m:
\uasn \otimes_\A \R \to \uan$ is surjective. By \cite{LTV}, Theorem 4.1 we have and section 1.8 we have
$$\mathrm{dim}\;\R_\K[\beta]=\#\{\m \in \mathcal{M}^{\text{per}}\;|\;\mathrm{dim}\;\m=\beta\}.$$
Moreover, it is well known (see Section 1.8) that 
\begin{align*}
&\mathrm{dim}\;\uksn[\alpha]=\#\{\m \in \mathcal{M}^{\text{ap}}\;|\;\mathrm{dim}\;\m=\alpha\},\\
&\mathrm{dim}\;\ukn[\gamma]=\#\{\m \in \mathcal{M}\;|\;\mathrm{dim}\;\m=\gamma\}.
\end{align*}
Thus,
$$\sum_{\alpha+\beta=\gamma}\mathrm{dim}\;\uksn[\alpha]+\mathrm{dim}\;\R_{\K}[\beta]=\mathrm{\dim}\;\ukn[\gamma].$$
This implies that $m$ is injective, and $ii)$ follows.\\
\hbox to1em{\hfill}To prove $iii)$, note that by definition, $\R=(\sum_i \f_i \uan)^\perp$. Let $x \in \R$, $y \in \sum_i \f_i\uan$ and $u \in \uan$. Then $yu, uy \in \sum_i \f_i\uan$ by $ii)$. Thus
\begin{alignat*}{2}
(\tilde{\Delta}(x),y\otimes u) &=(x,yu)&=0,\\
(\tilde{\Delta}(x),u \otimes y)&=(x, uy)&=0.
\end{alignat*}
Hence $\tilde{\Delta}(x) \in (\U^0\R)^{\otimes 2}$. Finally, it follows from $i)$ and section 1.5 that the map $\mathbf{U}^0\R \to \R, \; \mathbf{k}_\mathbf{a}u \mapsto u$ is an algebra homomorphism, which implies that $(\R,\Delta)$ is a bialgebra.
\qed

\paragraph{2.2} For $D \in \N$ set $\Gamma_D=\A[y_1,\ldots y_D]^{\Si_D}$, and let $\Gamma$ be the ring of symmetric functions in the variables $y_i, \;i\in \Z$. For $k \in \Z^*$, denote by $p_k^D, \;p_k$ the $k$-th power sum in $\Gamma_D$ and $\Gamma$. Recall that $\Gamma =\A[p_1,p_2,\ldots]$ and that $\Gamma$ is equipped with a canonical cocommutative bialgebra structure $\Delta_0: \Gamma \to \Gamma \otimes \Gamma, \; p_k \to p_k \otimes 1 + 1 \otimes p_k$ (see \cite{McD}, I, 5). Let us denote by $\rho_D: \uan \to \mathrm{End}\;\Lambda^D$ (resp. $\rho: \uan \to \mathrm{End}\;\La$) the representations of the Hall algebra on $\Lambda^D$ and $\Lambda^\infty$ (see Sections 1.10 and 1.11). Identify $\Gamma_D$ with $Z^-_D$ via $y_i \mapsto X_i^{-1}$ and let $\sigma_D: \Gamma_D \to \mathrm{End}\;\La$ (resp. $\sigma: \Gamma \to \mathrm{End}\;\La$) be the representations of the center of the affine Hecke algebras (see \cite{KMS} Section 1.).
\begin{prop} There exists a graded bialgebra isomorphism $i: (\R,\Delta) \stackrel{\sim}{\to} (\Gamma,\Delta_0)$ such that $\rho_{|\R}=\sigma \circ i$.
\end{prop}
\noindent \textit{Proof.} The action of $\uasn$ on $\La$ extends to a level 1 action of the whole quantum affine algebra $\U(\hat{\mathfrak{sl}}_n)$, which commutes to the action of $\Gamma$. Moreover, by \cite{KMS}, Proposition 2.3, the Fock space $\Lambda^\infty$ decomposes as
$$\La = L(\Lambda_0) \otimes_\A \Gamma$$
where $L(\Lambda_0)$ is the $\U(\hat{\mathfrak{sl}}_n)$-submodule generated by the highest weight vector $v_0$ in the irreducible representation of $\U_\K(\hat{\mathfrak{sl}}_n)$ with highest weight $\Lambda_0$.
\begin{lem} The actions of $\R$ and $\Gamma$ on $\La$ commute. 
\end{lem}
\noindent \textit{Proof.} Let $(a_D)_D$ be a family of operators such that $a_D \in \mathrm{End}\;\Lambda^D$. We say that $(a_D)_D$ satisfies property $(*)$ if there exists $N>0$ such that
\begin{equation*}
\begin{split}
\forall\; s,r \in \N,\;r\geq s+N,\;\; \forall\;i_1,\ldots, i_r,&\\
i_{l+1}=i_l-1 \;\mathrm{if}\;l\geq s \Rightarrow& \\
a_r(\x_{i_1} \wedge \ldots
\wedge \x_{i_r})=a_{s+N}(\x_{i_1}\wedge \ldots \wedge &\x_{i_{s+N}})\wedge
\x_{i_{s+N+1}}\wedge \ldots \wedge \x_{i_r}
\end{split}
\end{equation*}
It is easy to check that $(p_k^D)_D$ satisfies $(*)$ for any $k \in \N^*$ and that for any $u \in \uan$ the family $(\rho_D(u))_D$ satisfies $(*)$.
 As a consequence, the family $([\rho_D(u),p^D_k])_D$ satisfies $(*)$ with constant, say $N$. Given $\x_{i_1} \wedge \x_{i_2} \wedge \ldots$ with $i_l=1-l$ for $l \geq s$, we have
\begin{equation*}
[\rho(u),p_{k}](\x_\mathbf{i})=[\rho_{s+N}(u),p^{s+N}_{k}](x_{i_1}\wedge \ldots \wedge x_{i_{s+N}})\wedge x_{i_{s+N+1}} \wedge \ldots=0
\end{equation*}
since the actions of $\uan$ and $\Gamma_{s+N}$ on $\Lambda^{s+N}$ commute (see Section 1.10).\qed
\paragraph{}The proof of the following lemma will be given in the appendix (Section 4.1.)
\begin{lem}  The actions of $\R$ and $\U(\hat{\mathfrak{sl}}_n)$ on $\La$ commute.
\end{lem}
\paragraph{} It follows from Lemma 2.2 that for any $x \in \R$ there exists $\tilde{x} \in \Gamma$ such that $x.(v_0 \otimes 1)=v_0 \otimes \tilde{x}.$
 Let $i$ be the map $\R \to \Gamma, \;x \mapsto \tilde{x}$. Lemmas 2.1 and 2.2 imply that $\rho(x)=\sigma (i(x)) \in \mathrm{End}\;\La$. Moreover, by \cite{VV} Section 9.3 we have $\uan\vac=\La$. Hence $i$ is surjective. The map $i$ is a graded algebra homomorphism, where $\mathrm{deg}\;p_k=k\s$. Since $\R$ and $\Gamma$ have the same graded dimensions, it follows that $i$ is an isomorphism.\\
\hbox to1em{\hfill}Finally, we prove that $i$ is compatible with the bialgebra structures. For $D_1, D_2 \in \N$, consider the map 
$$\mu_{D_1,D_2}: \Lambda^{D_1} \otimes \Lambda^{D_2} \to \Lambda^{D_1+D_2},\qquad\wedge \x_\i \otimes \wedge \x_\mathbf{j}  \mapsto \wedge \x_{\i\mathbf{j}}$$
where $\i\mathbf{j}$ is the sequence $\i$ followed by the sequence $\mathbf{j}$. Then, by Section 1.11, (\ref{E:1103}) and the inclusion $\R \subset \bigoplus_{l} \U^-_n[l\s]$,
$$\mu_{D_1,D_2} \circ (\rho_{D_1} \otimes \rho_{D_2}) \circ \tilde{\Delta}=\mu_{D_1,D_2} \circ (\rho_{D_1} \otimes \rho_{D_2}) \circ \Delta=\rho_{D_1+D_2}.$$
Let $k \in \N$ and put $c_k=i^{-1}(p_k)$. Choose $N \gg 0$ such that property $(*)$ is satisfied by $(\rho_D(u))_D$ for any $u \in \bigoplus_{l \leq k} \R[l\tau]$ with constant $N$. Thus,
\begin{equation*}
\begin{split}
\forall\; s,r,D \in \N,\;r\geq s+D, \forall\;\i \in \Z^D\mathrm{\;such\;that\;}
i_{l+1}=i_l-1 \;\mathrm{if}\;l&\geq s,\\
 c_k(\x_{i_1}\wedge \ldots
\wedge \x_{i_D})=p_k^{s+D}(\x_{i_1}\wedge \ldots \wedge \x_{i_{s+D}})\wedge
\x_{i_{s+D+1}}&\wedge \ldots \wedge \x_{i_r}.
\end{split}
\end{equation*}
In particular, if $\i \in \Z^{D_1}$ and $\mathbf{j} \in \Z^{D_2}$ both satisfy the condition in $(*)$ then
$$\mu_{D_1,D_2}(\rho_{D_1} \otimes \rho_{D_2}) \Delta(c_k)(\wedge \x_\i \otimes \wedge \x_\mathbf{j})=\mu_{D_1,D_2}(1 \otimes c_k + c_k \otimes 1) (\wedge \x_\i \otimes \wedge \x_\mathbf{j})$$
It is easy to see that this implies that $\Delta(c_k)=1 \otimes c_k + c_k \otimes 1$ as desired.
\qed
 \paragraph{}Propositions 2.1 and 2.2 together imply
\begin{theo} There exists a graded algebra isomorphism $\uan \simeq \uasn \otimes_\A \R$ where $\R\simeq \A[\mathbf{z}_1,\mathbf{z}_2\ldots ]$ with $\mathrm{deg}\;\mathbf{z}_i=i\s \in \N^I$.
\end{theo}

\paragraph{2.3} Let us denote by $\rho'_D$ and $\sigma'_D$ the representations of $\uan$ and $\Gamma_D$ on $\TD$. Then in fact:
\begin{prop} There exists graded algebra morphisms $i_D: \R \twoheadrightarrow \Gamma_D$, $D \in \N^*$ such that $i=\underset{\longleftarrow}{\mathrm{lim}}\; i_D$ and $\rho'_D=\sigma'_D \circ i_D$.
\end{prop}
\noindent
\textit{Proof.} It follows from (\ref{E:1102}) and from the fact that $\R$ belongs to the center of $\uan$ that, for any $k \in \N$ there exists $k'\in \N$ such that
$$\forall\;i \in \Z,\;\; c_k.\x_i=\x_{i+nk'}.$$
Then, by Prop. 2.2, for all $D \in \N$,
$$\forall\;\i \in \Z^D,\;\; \rho'_D(c_k) (\x_{i_1} \otimes \ldots \otimes \x_{i_D}) = \sigma'_D(p_{k'}^D)(\x_{i_1} \otimes \ldots \otimes \x_{i_D}).$$
In particular, $\rho_D(c_k)=\sigma_D(p_{k'}^D)$ for all $D \in \N$. This implies $k=k'$ and proves the proposition. \qed
\paragraph{2.4} Let $(\b_\m^*)_{\m \in \mathcal{M}}$ be the dual basis of $(\b_\m)_{\m \in \mathcal{M}}$ with respect to $(\,,\,)$.
\begin{prop} We have $\R=\bigoplus_{\m \in \mathcal{M}^{\mathrm{per}}} \C[v,v^{-1}] \b^*_{\m}$.
\end{prop}
\noindent
\textit{Proof.} By definition we have 
$$\R=\bigcap_i \mathrm{Ker}\;e'_i=( \sum_i \f_i \U^-_{\k,n})^\perp.$$
From \cite{Lb}, Th.14.3.2 and from the geometric description of $\U^-_n$ in terms of Frobenius traces of perverse sheaves on $E_{V_\d}$, $\d \in \N^{\Z/n\Z}$ (see \cite{VV}) it follows that
 $$\sum_i \f_i \U^-_{n}=\bigoplus_{i,\m \in \tilde{f}_i \mathcal{C}} \C[v,v^{-1}]\b_\m=\bigoplus_{\m \not\in \mathcal{M}^\mathrm{per}} \C[v,v^{-1}] \b_\m.$$
This proves the Proposition. \qed
\noindent
\section{Proof of the Varagnolo-Vasserot conjecture}
\paragraph{3.1} Set $\mathcal{L}_{\Lambda}=\bigoplus_\lambda \S |\lambda \rangle \subset \Lambda^\infty$. Leclerc and Thibon have defined a semilinear involution $a \mapsto \overline{a}$ on $\Lambda^\infty$ such that
\begin{enumerate}
\item[i)]$\overline{\vac}=\vac,$
\item[ii)]$\overline{ua}=\overline{u}\;\overline{a}\;\mathrm{for\;all}\; u \in \uan,\; a \in \Lambda^\infty,$
\item[iii)]$\overline{p_k.a}=p_k.\overline{a}\;\mathrm{for\;all}\; k \in \N^*,\; a \in \Lambda^\infty,$
\end{enumerate}
(see \cite{LT}, \cite{VV}).\\
\hbox to1em{\hfill}For $\lambda \in \Pi$, set $\m(\lambda)=\sum_i [1-i,\lambda_i-i] \in \mathcal{M}$. To simplify notations, put $\f_\lambda=\f_{\m(\lambda)}$ and $\b_\lambda=\b_{\m(\lambda)}$. If $\lambda \in \Pi$, let $n\lambda$ be the partition $((\lambda_1)^n,(\lambda_2)^n,\ldots)$. Thus $\mathcal{M}^{\mathrm{per}}=\{\m(n\lambda),\;\lambda \in \Pi\}$. Leclerc and Thibon introduced in \cite{LT} two canonical bases $\mathbf{B}^{\pm}=\{\b^{\pm}_\lambda,\;\lambda \in \Pi\}$ of $\La$ caracterized by
\begin{equation}\label{E:3}
 \overline{\b^{\pm}_\lambda}=\b^{\pm}_\lambda,\;\qquad \b^+_\lambda \in |\lambda\rangle + v\bigoplus_{\mu <\lambda} \S |\mu\rangle,\qquad \b^-_\lambda \in |\lambda\rangle + v^{-1}\bigoplus_{\mu <\lambda} \overline§ |\mu\rangle.
\end{equation}
The following was conjectured in \cite{VV} and is the main result of
this paper :
\begin{theo} For all $\lambda \in \Pi$ we have $\b_\lambda \vac=\b^+_\lambda$.
\end{theo}
\noindent The rest of this section is devoted to the proof of this theorem.

\paragraph{3.2} Recall that a partition $\lambda=(\lambda_1 \geq \lambda_2 \geq \ldots)$ is called $n$-regular if $\lambda_{i}>\lambda_{n+i}$ for all $i$ such that $\lambda_i \neq 0$. Let $\Pi^{\mathrm{reg}}$ be the set of all $n$-regular partitions. We first show
\begin{lem} We have $\b_{\lambda}|0\rangle=\b^+_\lambda$ if $\lambda\in \Pi^{\mathrm{reg}}$.
\end{lem}
\noindent
\textit{Proof.} Consider the scalar product $(\;,\;)$ on $\La$ for which $\{|\lambda\rangle\}$ is an orthonormal basis. Recall that $\La$ is isomorphic to $L(\Lambda_0) \otimes_\A \Gamma$ as a $\uasn \otimes_\A \Gamma$-module. It is shown in \cite{LT2} that the restriction of $(\,,\,)$ and of the involution $a \mapsto \overline{a}$ to $L(\Lambda_0)$ coincide with the Kashiwara scalar product and involution defined on any simple integrable $\uan$-module (\cite{K}, Sections 2 and 6). Thus the lower crystal basis of $L(\Lambda_0)$ is a subset of $\{\pm \b^+_{\lambda}\}$. Note that $\m(\lambda) \in \mathcal{M}^{\text{ap}}$ if and only if $\lambda$ is $n$-regular. Therefore, by Section 1.8 and the general theory of canonical bases
\begin{equation}\label{E:321}
\lambda \in \Pi^{\mathrm{reg}} \Rightarrow \b_{\lambda}|0\rangle \subset \{\pm \b^+_{\lambda}\}.
\end{equation}
 Moreover, by \cite{VV}, Section 9.2, for any $\lambda \in \Pi$ and any orbit $O \subset \overline{O}_\lambda \backslash O_\lambda$ we have
\begin{equation}\label{E:322}
\f_\lambda \vac\in |\lambda \rangle + \bigoplus_{\mu<\lambda} \A |\mu\rangle, \qquad \f_O \vac\in \bigoplus_{\mu<\lambda} \A |\mu\rangle.
\end{equation}
It is now clear from (\ref{E:3}) and (\ref{E:321}), (\ref{E:322}) that $\b_{\lambda}|0\rangle=\b^+_\lambda$ if $\lambda$ is an $n$-regular partition.\qed

\paragraph{3.3} Set $\mathcal{L}^{\mathrm{reg}}=\bigoplus_{\m \in \mathcal{M}^{\mathrm{ap}}} \S \b_\m$. It is known that $\mathcal{L}^{\mathrm{reg}}$ is the smallest $\S$-submodule of $\uan$ containing $1$ and stable by the operators $\tilde{f}_i$, $i \in I$ (c.f. \cite{K} and Section 1.8). Set $\mathcal{L}_\R=\mathcal{L} \cap \R$. If $V$ is any $\S$-module $V$, we let $\overline{V}=V \otimes_\S \C[[v]]$ be its completion with respect to the $v$-adic topology.
\begin{lem} The multiplication defines a graded isomorphism
$$\overline{\mathcal{L}^{\mathrm{reg}}}\otimes_{\C[[v]]} \overline{\mathcal{L}}_\R \stackrel{\sim}{\to} \overline{\mathcal{L}}.$$
\end{lem}
\noindent
\textit{Proof.} Every multisegment $\mathbf{n}$ decomposes in a unique
way as $\mathbf{n}=\mathbf{p}+\mathbf{a}$ where $\mathbf{p}\in
\mathcal{M}^{\mathrm{per}}$ and $\mathbf{a} \in
\mathcal{M}^{\mathrm{ap}}$. Since $\text{b}_{\mathbf{a}}$ belongs to
the connected component of the crystal graph $\mathcal{C}$ containing
$1$, there exists a sequence $i_1,\ldots i_r$ such that
$\tilde{f}_{i_1} \cdots \tilde{f}_{i_r} 1 \equiv
\mathbf{f}_{\mathbf{a}}\;(\mathrm{mod}\;v\mathcal{L})$. Moreover, by
Section 1.7 there exists $x \in \mathcal{L}_\R$ such that $x \equiv
\f_\mathbf{p}\;(\mathrm{mod}\;v\mathcal{L})$. Since the left
multiplication by $x$ commutes with the $\tilde{f}_i$ (see the proof of Proposition 2.1), we have
\begin{alignat*}{2}
x \tilde{f}_{i_1} \cdots \tilde{f}_{i_r}.1 &=\tilde{f}_{i_1} \cdots \tilde{f}_{i_r}\cdot x\qquad&\qquad &\\
&\equiv \tilde{f}_{i_1} \cdots \tilde{f}_{i_r}\cdot \mathbf{f}_{\mathbf{p}} & &(\mathrm{mod}\;v\mathcal{L})\\
&\equiv \mathbf{f}_{\mathbf{n}} &  &(\mathrm{mod}\;v\mathcal{L})
\end{alignat*}
Hence the multiplication map induces an isomorphism modulo $v$, and hence an isomorphism over $\C[[v]]$.\qed
\paragraph{3.4} For $\nu \in \Pi$ let $s_\nu \in \Gamma$ be the Schur polynomial and set $a_\nu=i^{-1}(s_\nu) \in \R$. Then $\{a_\nu,\;\nu \in \Pi\}$ is an $\A$-basis of $\R$.
\begin{lem}The following holds :
\begin{enumerate}
\item[i)] $v^{(n-1)|\mu|}a_\mu \mathcal{L}_\Lambda \subset \mathcal{L}_\Lambda$,
\item[ii)] $a_\mu \vac \in (-v)^{(n-1)|\mu|}\big( |n\mu\rangle + v \mathcal{L}_\Lambda\big)$.
\item[iii)] More generally, let $\lambda \in \Pi$ and write $\m(\lambda)=\m(\lambda')+\m(n\mu)$ where $\lambda'\in \Pi^{\mathrm{reg}}$. Then
$$a_\mu |\lambda'\rangle \in (-v)^{-(n-1)|\mu|}\big(|\lambda \rangle + v \mathcal{L}_\Lambda\big).$$
\end{enumerate}
\end{lem}
\noindent
\textit{Proof.} Statement $i)$ follows from \cite{LT2}, Theorem 6.3. Statements $ii)$ and $iii)$ are proved as in \cite{LT2}, Theorem 6.7. \qed

\paragraph{}Let us denote by $<$ the order on multisegments such that $\m \leq \mathbf{n}$ if $O_{\mathbf{m}} \subset \overline{O_\mathbf{n}}$.

\begin{prop} There holds
\begin{enumerate}
\item[i)] $a_\lambda \in (-v)^{-(n-1)|\lambda|}(\f_{n\lambda} +
  v\mathcal{L}),$
\item[ii)] $\mathcal{L}_\R\cdot \mathcal{L}_\Lambda \subset \mathcal{L}_\Lambda$ and $\mathcal{L} \vac \subset \mathcal{L}_\Lambda$.
\end{enumerate}
\end{prop}
\noindent
\textit{Proof.} Let $i)_k$ be statement $i)$ restricted to all $\lambda
\in \Pi $ with $|\lambda| \leq k$, and let $ii)_k$ be statement $ii)$
restricted to $\bigoplus_{k'\leq k}\mathcal{L}_\R[k'\tau]$, and $\mathcal{L}[\mathbf{d}]$, $\mathbf{d}=(d_1,\ldots d_n),\; d_i \leq k$. We will prove $i)_k$ and $ii)_k$ by induction. The case $k=1$ is a consequence of Lemma 3.2 and the following formula :
\begin{equation}\label{E:341}
a_{(1)}=\sum_{\underset{\mathbf{dim}\;\mathbf{m}=\s}{\mathbf{m}=[a_1;l_1)+\ldots+[a_r;l_r)}}(v-v^{-1})^{(n-1)-\sum_{i}(l_i-1)}\mathbf{f}_{\mathbf{m}}.
\end{equation}
Indeed, let $x$ denote the r.h.s of (\ref{E:341}). A direct computation using \cite{LTV}, Proposition 4.1, shows that $x \in \R$. Since $\mathrm{dim}\;\R[\tau]=1$, we have
$a_{(1)}=cx$ for some $c \in \K$. Using \cite{LTV} Theorem 6.3, we see
that the coefficient of $|(n)\rangle$ in $a_{(1)}\vac$ is equal to
$1$. On the other hand, by (\ref{E:322}), we have
\begin{align*}
cx \vac &=c\sum_{\underset{\mathbf{dim}\;\mathbf{m}=\s}{\mathbf{m}=[a_1;l_1)+\ldots+[a_r;l_r)}} (v-v^{-1})^{(n-1)-\sum_i (l_i-1)} \f_\m \vac\\
&\in c \f_{(n)}\vac + \bigoplus_{\mu <(n)} \A |\mu \rangle\\
&\in c |(n)\rangle +  \bigoplus_{\mu <(n)} \A |\mu \rangle.
\end{align*}
Therefore $c=1$ and (\ref{E:341}) is proved.

\paragraph{}For $\lambda,\mu \in \Pi$ and $x \in \R$ let $\Delta_{\lambda,\mu}(x) \in \A$ be the coefficient of $\f_{n\lambda} \otimes
\f_{n\mu}$ in $\Delta(x)$, where $\Delta(x)$ is expressed in the basis $(\f_\m \otimes \f_{\m'})_{\m,\m'}$. For $\lambda, \mu,\nu \in \Pi$, $|\lambda|=|\mu|+|\nu|$ we let $c^\lambda_{\mu,\nu} \in \N$ be the Littlewood-Richardson multiplicity (see \cite{McD}, Section 5).
\begin{lem} For all $\lambda,\mu \in \Pi$, $l \in \N^*$ and $\m \in
  \mathcal{M}\backslash \mathcal{M}^{\mathrm{per}}$ we have
$$\Delta_{(1)^l,\mu}(\f_{n\lambda})\in c^\lambda_{(1)^l,\mu} + v \S,\qquad \mathrm{and} \qquad
\Delta_{(1)^l,\mu}(\f_{\m})\in v \S.$$
\end{lem}
\noindent
\textit{Proof.} See the appendix. \qed
\paragraph{}Now let $k>1$ and suppose that $i)_{k-1}$ and $ii)_{k-1}$ hold. Let $d_\nu \in \Z$, $\nu \in \Pi$, be such that $\{v^{d_\nu}a_\nu\}$ is a $\S$-basis of $\mathcal{L}_\R$. It follows from the crystal graph $\mathcal{C}$ of $\U^-_n$ that $e'_i(\text{b}_\m)=0$ for all $i \in I$ if and only if $\m \in \mathcal{M}^{\mathrm{per}}$. Hence
\begin{equation}\label{E:n31}
a_\nu=(-v)^{-d_\nu}\left(\sum_{\sigma\in \Pi} \alpha_\sigma \f_{n\sigma}
  +\sum_{\mathbf{l}} \beta_{\mathbf{l}} \f_{\mathbf{l}\in \mathcal{M}}\right)
\end{equation}
for some $\alpha_\sigma \in \C$ and $\beta_{\mathbf{l}} \in v \S$. Thus, by Lemma 3.4, 
\begin{equation}\label{E:343}
\Delta_{(1)^l,\mu}(a_\nu) \in (-v)^{-d_\nu}(\sum_\sigma \alpha_\sigma c^\sigma_{(1)^l,\mu} + v\S).
\end{equation}
On the other hand, by \cite{McD}, Section 5.3 and Proposition 2.3, we have
\begin{equation}\label{E:frt}
\Delta(a_\nu)=\sum_{|\lambda|+|\mu|=k} c^\nu_{\lambda,\mu} a_\lambda \otimes a_\mu.
\end{equation}
Using the induction hypothesis $i)_{k-1}$ we obtain, for $|\lambda|,|\mu|\neq 0$
\begin{equation}\label{E:342}
\Delta_{\lambda,\mu}(a_\nu) \in (-v)^{-(n-1)|\nu|}(c^\nu_{\lambda,\mu} + v\S).
\end{equation}
For any $\nu \in \Pi$ there exists some $l \in \N^*$ and $\mu \in \Pi$ such that $c^\nu_{(1)^l,\mu}=1$ (see \cite{McD}, (5.17)). Combining (\ref{E:343}), (\ref{E:342}), we see that $d_\nu \geq (n-1)|\nu|$. In particular, it follows from Lemma 3.3 $i)$ that $v^{d_\nu} a_\nu . \mathcal{L}_\Lambda \subset \mathcal{L}_\Lambda$. Since $\{v^{d_\nu}
 a_\nu\}_{|\nu|=k}$ is an $\S$-basis of $\mathcal{L}_\R[k\tau]$ we obtain the first statement of $ii)_{k}$. The second part of Statement $ii)_k$ follows from Lemmas 3.1 and 3.2. Then $ii)_k$ implies
\begin{equation*}
\begin{split}
a_\nu \vac \in& (-v)^{-d_\nu} \left( \sum_\sigma \alpha_\sigma \f_{n\sigma} + v \mathcal{L}[k\tau] \right) \vac\\
\in& (-v)^{-d_\nu}\left( \sum_\sigma \alpha_\sigma \f_{n\sigma} \vac + v\mathcal{L}_\Lambda\right).
\end{split}
\end{equation*}

Using (\ref{E:322}), $ii)_k$ and Lemma 3.3 $ii)$, we get 
$$d_\nu=(n-1)|\nu|,\qquad \alpha_{\nu}=1,\qquad \alpha_{\sigma} \neq 0 \Rightarrow \sigma \leq \nu.$$
Finally, combining (\ref{E:342}) and (\ref{E:343}) now yields
$$\forall\;l\in \N,\;\forall\; \mu \in \Pi,\qquad \sum_{\sigma < \nu} c^\sigma_{(1)^l,\mu} \alpha_\sigma=0.$$
For any fixed $\nu$, this system is nondegenerated and admits the
unique solution $(\alpha_\sigma)_{\sigma < \nu}=0$. Indeed let
$\langle\;,\;\rangle$ be the nondegenerate symmetric bilinear form on
$\Gamma$ for which $\{s_\lambda\}$ is an orthonormal basis, and set
$X=\sum_\sigma \alpha_\sigma s_\sigma$. Then for any $\mu$ and $l$ we
have
$$\langle s_{(1)^l}s_\mu,X\rangle=\langle \sum_\sigma c^\sigma_{(1)^l,\mu}s_\sigma,X\rangle=\sum_{\sigma}
c^\sigma_{(1)^l,\mu}\alpha_\sigma=0.$$
Moreover $\langle s_{(k)},X\rangle=0$ since $X \in
\bigoplus_{\sigma <\nu} \A s_\sigma \subset \bigoplus_{\sigma<(k)} \A s_\sigma$. But
$$\Gamma[k]=\A s_{(k)} \oplus \sum_{l=1}^{k-1} s_{(1)^l}\Gamma[k-l].$$
Hence $X=0$. Statement $i)_k$ is proved and the induction is complete.\qed

\paragraph{3.5}\textit{Proof of Theorem 3.1.} Let $\lambda \in
 \Pi$. The multisegment $\m(\lambda)$ decomposes in a unique way as
 $\m(\lambda)=\m(\lambda')+\m(n\mu)$ for some partitions $\lambda' \in
 \Pi^{\mathrm{reg}}$ and $\mu \in \Pi$. Moreover, from Section 1.7 the element  $\text{b}_\lambda$ is in the connected component of the crystal
 graph $\mathcal{C}$ containing $\text{b}_{n\mu}$. Hence, by Proposition 3.4 $i)$ and the proof of Lemma 3.2, we have
$$\b_\lambda \equiv \f_\lambda \equiv (-v)^{(n-1)|\mu|}a_\mu \f_{\lambda'}\;(\mathrm{mod}\;v\mathcal{L}).$$
Then, by Proposition 3.4 $ii)$, Lemma 3.3 $i)$ and the fact that $\b_{\lambda'}-\f_{\lambda'} \in v \mathcal{L}$, we get 
$$(-v)^{(n-1)|\mu|}a_\mu \f_{\lambda'} \vac \equiv (-v)^{(n-1)|\mu|}a_\mu \b_{\lambda'} \vac \;(\mathrm{mod}\;v\mathcal{L}_\Lambda). $$
Finally, by Lemma 3.1, and Lemma 3.3 $i)$ and $iii)$, we have 
\begin{align*}
 (-v)^{(n-1)|\mu|}a_\mu \b_{\lambda'} \vac &\equiv  (-v)^{(n-1)|\mu|}a_\mu \b^+_{\lambda'} \vac\;(\mathrm{mod}\;v\mathcal{L}_\Lambda) \\
&\equiv (-v)^{(n-1)|\mu|}a_\mu |\lambda'\rangle \;(\mathrm{mod}\;v\mathcal{L}_\Lambda)\\
&\equiv |\lambda \rangle \;(\mathrm{mod}\;v\mathcal{L}_\Lambda).
\end{align*}
 Thus $\b_\lambda \vac \in |\lambda\rangle + v\bigoplus_{\mu <\lambda}
\S |\mu\rangle$. Finally,
$\overline{\b_\lambda\vac}=\overline{\b_\lambda}\,\overline{\vac}=\b_\lambda
\vac$ by Section 3.1. Hence $\b_\lambda \vac=\b^+_\lambda$ as desired. \qed
\noindent
\section{An analogue of the Varagnolo-Vasserot conjecture for higher-level Fock spaces}
\paragraph{}In this section we sketch the generalization of Theorem 3.1 to the case of the higher-level Fock spaces. We use
the definitions and notations of \cite{U}.
\paragraph{5.1} Let $l>1$ and $s \in \Z$. Let $\Lambda^{s+\infty}$ be the semi-infinite wedge product of levels $l$ and $n$
and charge $s$ (see \cite{U}, Section 4.1). Let 
$$\Z^l(s)=\{\mathbf{s}_l=(s_1,\ldots s_l) \in \Z^l\;|\sum_i s_i=s\}.$$
Recall that $\Lambda^{s+\infty}$ is equipped with a distinguished $\A$-basis $\{|\lambda_l,\s_l\rangle\}$ where
$\lambda_l=(\lambda_1,\ldots \lambda_l) \in \Pi^l$ and $\s_l \in \Z^l(s)$. It is endowed with three commuting left actions:
$\rho_{l,s}: \U(\hat{\mathfrak{sl}}_n) \to \mathrm{End}\;(\Lambda^{s+\infty})$, $\rho'_{l,s}: \U'(\hat{\mathfrak{sl}}_l) \to
\mathrm{End}\;(\Lambda^{s+\infty})$ where $\U'(\hat{\mathfrak{sl}}_l)$ denotes the Lusztig integral form of the quantum
affine algebra of type $A^{(1)}_{l-1}$ with quantum parameter $v^{-1}$ and the action of a Heisenberg algebra $\mathcal{H}$
generated by operators $B_m$, $m \in \Z^*$ (see \cite{U}, Sections 4.2 and 4.3). Moreover, $\Lambda^{s+\infty}$ is an
integrable module for $\U(\hat{\mathfrak{sl}}_n)$ and $\U'(\hat{\mathfrak{sl}}_l)$. We denote by $\mathcal{H}^-$ (resp.
$\mathcal{H}^+$) the subalgebra of $\mathcal{H}$ generated by $B_{-m}$, $m \in \N$ (resp. $B_m$, $m \in \N$).\\
\hbox to1em{\hfill}Set $\Lambda^{s+\infty}_\K=\Lambda^{s+\infty}\otimes_{\A} \K$. The Fock space
$\Lambda^{s+\infty}_\K$ decomposes under these actions as follows (see \cite{TU}, Theorem 4.10). Set
$$\mathcal{A}^n_l=\{\s_l
\in
\Z^l(s)\;|s_1 \geq s_2 \geq
\ldots
\geq s_l, s_1-s_l \leq n\}.$$
Then for $\s_l \in \mathcal{A}^n_l$ the vector $|0,\s_l\rangle$ is singular for
$\U(\hat{\mathfrak{sl}}_n)$, $\U'(\hat{\mathfrak{sl}}_l)$ and $\mathcal{H}$, i.e we have
$$ \U(\hat{\mathfrak{sl}}_n)^+|0,\s_l\rangle=\U'(\hat{\mathfrak{sl}}_l)^+|0,\s_l\rangle=\mathcal{H}^+|0,\s_l\rangle=0,$$
and
$$\Lambda_\K^{s+\infty}=\bigoplus_{\s_l \in
\mathcal{A}^n_l}\U(\hat{\mathfrak{sl}}_n)\cdot\mathcal{H}\cdot\U'(\hat{\mathfrak{sl}}_l)\;|0,\s_l\rangle.$$ Moreover, for
any
$\s_l
\in
\Z^l(s)$, $\Lambda_{\s_l}=\bigoplus_{\lambda_l \in \Pi^l} \A |\lambda_l,\s_l\rangle$ is a
$\U(\hat{\mathfrak{sl}}_n)\cdot\mathcal{H}$-module, and $|0,\s_l\rangle$ generates an irreducible integrable
$\U(\hat{\mathfrak{sl}}_n)$-module $V_{\s_l}$ of highest weight $\Lambda=\Lambda_{s_1}+\Lambda_{s_2}+\ldots +
\Lambda_{s_l}$. Here $\Lambda_i$, $0 \leq i \le n-1$ is the $ith$ fundamental weight of $\hat{\mathfrak{sl}}_n$ and
we set $\Lambda_i=\Lambda_j$ if $i \equiv j\;(\mathrm{mod}\;n)$.
\paragraph{}Set $\mathcal{L}_{\Lambda^{s+\infty}}=\bigoplus_{\lambda_l,\s_l} \S |\lambda_l,\s_l\rangle$ and
$\mathcal{L}_{\s_l}=\bigoplus_{\lambda_l} \S |\lambda_l, \s_l\rangle$ for any $\s_l \in \Z^l(s)$. Let
$\text{b}_{\lambda_l,\s_l}$ denote the image of $|\lambda_l,\s_l\rangle$ in
$\mathcal{L}_{\Lambda^{s+\infty}}/v\mathcal{L}_{\Lambda^{s+\infty}}$. Set
$\text{B}_{\Lambda^{s+\infty}}=\{\text{b}_{\lambda_l,\s_l}\;|\;\lambda_l \in \Pi^l,\;\s_l \in \Z^l(s)\}$ and for any $\s_l
\in \Z^l(s)$ set $\text{B}_{\s_l}=\{\text{b}_{\lambda_l,\s_l}\;|\;\lambda_l \in \Pi^l\}$ .
The following is proved in
\cite{JMMO}.
\begin{prop} The couple $(\mathcal{L}_{\Lambda^{s+\infty}},\text{B}_{\Lambda^{s+\infty}})$ is a crystal basis of the
$\U(\hat{\mathfrak{sl}}_n)$-module $\Lambda^{s+\infty}$. 
\end{prop}
\noindent
The crystal graph structure of $\text{B}_{\Lambda^{s+\infty}}$ is explicitely described in \cite{JMMO}. 
\paragraph{}In \cite{U}, Uglov has defined a semilinear involution $a \mapsto \overline{a}$ on $\Lambda^{s+\infty}$
satisfying
\begin{enumerate}
\item[i)]$\overline{|0,\s_l\rangle}=|0,\s_l\rangle$ for any $\s_l \in \Z^l(s)$,
\item[ii)]$\overline{ua}=\overline{u}\;\overline{a}\;\mathrm{for\;all}\; u \in \U(\hat{\mathfrak{sl}}_n),\U'(\hat{\mathfrak{sl}}_l),\; a \in \Lambda^{s+\infty}$,
\item[iii)]$\overline{B_{-m}a}=B_{-m}\overline{a}\;\mathrm{for\;all}\; m \in \N^*,\; a \in \Lambda^{s+\infty}$.
\end{enumerate}
Uglov also introduced two canonical bases $\{\b^{\pm}_{\lambda_l,\s_l}\}_{\lambda_l \in \Pi^l,\s_l \in \Z^l(s)}$
caracterized by the following properties
$$\overline{\b^{\pm}}_{\lambda_l,\s_l} =\b^{\pm}_{\lambda_l,\s_l},$$
$$\b^+_{\lambda_l,\s_l} \in |\lambda_l,\s_l\rangle +v\bigoplus_{\mu_l,\mathbf{t}_l}\S |\mu_l,\mathbf{t}_l\rangle ,\qquad \b^-_{\lambda_l,\s_l} \in |\lambda_l,\s_l\rangle +v^{-1}\bigoplus_{\mu_l,\mathbf{t}_l} \overline{\S}|\mu_l,\mathbf{t}_l\rangle.$$
The set $\mathbf{B}^+_{\s_l}=\{\b^+_{\lambda,\s_l}\;|\lambda_l \in \Pi^l\}$ is a basis of $\Lambda_{\s_l}$ which contains the lower canonical basis of the irreducible $\U(\hat{\mathfrak{sl}}_n)$-module $V_{\s_l}$.
\paragraph{5.2} By Theorem 2.2, we can extend the action of $\U(\hat{\mathfrak{sl}}_n)$ on $\Lambda^{s+\infty}$ to an action of $\uan$ by setting $\rho_{l,s}(i^{-1}(p_k))=B_{-k}$. Recall that $\mathcal{L}=\bigoplus_{\m} \S \f_\m$ and $\mathcal{L}_\R=\mathcal{L} \cap \R$ denote the integral lattices in $\uan$ and $\R$ respectively, and that $a_\mu=i^{-1}(s_\mu) \in \R$ denotes the Schur polynomial associated to $\mu \in \Pi$ (see Section 3.4).
\begin{prop} We have, for any $\mu \in \Pi$ and $\s_l \in \Z^l(s)$,
\begin{enumerate}
\item[i)] $(-v)^{(n-1)|\mu|}a_\mu \mathcal{L}_{\Lambda^{s+\infty}} \subset
\mathcal{L}_{\Lambda^{s+\infty}}$ and $(-v)^{(n-1)|\mu|}a_\mu \text{B}_{\s_l}\subset \text{B}_{\s_l}$,
\item[ii)] $\mathcal{L}|0,\s_l\rangle \subset \mathcal{L}_{\Lambda^{s+\infty}}$.
\end{enumerate}
\end{prop}
\noindent
\textit{Proof.} It will be convenient to use the dual indexation of elements of the basis $\{|\lambda_l,\s_l\rangle\}$ by pairs
$(\lambda_n,\s_n)$ where $\lambda_n \in \Pi^n$ and 
$$\s_n\in \Z^n(s)=\{(s_1,\ldots s_n) \in \Z^n\;|\;\sum_i s_i=s\}$$ as explained
in \cite{U}, Section 4.1. In particular we set
$\Lambda_{\s_n}=\bigoplus_{\mu_n}\A|\mu_n,\s_n\rangle$ and $\mathcal{L}_{\Lambda_{\s_n}}=\bigoplus_{\mu_n}
\S |\mu_n,\s_n\rangle$. Now let $\mu_n \in \Pi^n$, $\s_n \in \Z^n(s)$. It follows from \cite{U}, Corollary 5.6, i') that
$(-v)^{(n-1)|\lambda|}a_\lambda |\mu_n,\s_n\rangle \in \mathcal{L}_{\Lambda_{\s_n}}$ if $(\mu_n,s_n)$ is
$l|\lambda|$-dominant (see \cite{U}, Section 5.1). Let $\tilde{e}_i,\tilde{f}_i$, $i=0,\ldots n-1$ be the Kashiwara operators
corresponding to the $\U(\hat{\mathfrak{sl}}_n)$-action on $\Lambda^{s+\infty}$. By \cite{U}, Corollary 4.9, there exists a
sequence $j_1,\ldots j_r$ and operators $\tilde{x}_{1},\ldots \tilde{x}_{r}$ with $\tilde{x}_i \in
\{\tilde{e}_{j_i},\tilde{f}_{j_i}\}$ such that $\tilde{x}_{1}\ldots \tilde{x}_r |\mu_n,\s_n\rangle$ is a sum of
$l|\lambda|$-dominant vectors. Then
$$\tilde{x}_1\ldots \tilde{x}_r(-v)^{(n-1)|\lambda|}a_\lambda |\mu_n,\s_n\rangle=(-v)^{(n-1)|\lambda|}a_\lambda
\tilde{x}_1\ldots
\tilde{x}_r|\mu_n,\s_n\rangle \in \mathcal{L}_{\Lambda^{s+\infty}}.$$
But $\tilde{x}_1\ldots \tilde{x}_r$ defines an isomorphism $\Lambda_{\s_n} \to \tilde{x}_1\ldots
\tilde{x}_r\Lambda_{\s_n}$ which restricts to an isomorphism
$$\mathcal{L}_{\Lambda_{\s_n}} \otimes_\A \C[[v]] \to \tilde{x}_1\ldots
\tilde{x}_r \mathcal{L}_{\Lambda_{\s_n}} \otimes_\A \C[[v]]=(\mathcal{L}_{\Lambda^{s+\infty}}\cap  \tilde{x}_1\ldots
\tilde{x}_r \Lambda_{\s_n}) \otimes_\A \C[[v]].$$
It follows that $(-v)^{(n-1)|\lambda|}a_\lambda  |\mu_n,\s_n\rangle \in \mathcal{L}_{\Lambda_{\s_n}}$. Hence
$(-v)^{(n-1)|\lambda|}a_\lambda\mathcal{L}_{\Lambda^{s+\infty}}\subset
\mathcal{L}_{\Lambda^{s+\infty}}$ and the first statement of $i)$ is proved. The second statement  of $i)$ is proved in the
same way using Lemma 3.3 and \cite{U} Corollary 5.6 i').\\
\hbox to1em{\hfill}By the general theory of canonical bases,
$\b_\m^+|0,\s_l\rangle
\in
\mathbf{B}^+_{\s_l}
\cap V_{\s_l}$ for any $\m \in \mathcal{M}^{\mathrm{ap}}$ (see \cite{U}, Section 4.4). Hence $\mathcal{L}^{\mathrm{reg}}
|0,\s_l\rangle \subset
\mathcal{L}_{\Lambda_{\s_l}}$. Statement $ii)$ is now a consequence of Lemma 3.2.\qed

\begin{theo}For any $\s_l \in \Z^l(s)$ we have $\mathbf{B}|0,\s_l\rangle \subset \mathbf{B}^+_{\s_l}\cup \{0\}$.\end{theo}
\noindent
\textit{Proof.} By \cite{U}, Section 4.4, the basis $\mathbf{B}^+_{\s_l}$ contains the upper canonical basis of the irreducible
$\U(\hat{\mathfrak{sl}}_n)$-module $V_{\s_l}$. Thus, by Section 1.11
$$\b_\m |0,\s_l\rangle \in \mathbf{B}^+_{\s_l}\cup\{0\},\qquad \forall\;\m \in \mathcal{M}^{\mathrm{ap}}.$$
Now any $\m \in \mathcal{M}$ decomposes as $\m=\m(\lambda')+\m(n\mu)$ for some $\lambda' \in \Pi^{\mathrm{reg}}$ and
$\mu \in \Pi$. By Lemma 3.3 and Proposition 3.4,
$$\b_\m \equiv (-v)^{(n-1)|\mu|}a_\mu \f_{\lambda'}\;(\mathrm{mod}\;v\mathcal{L}).$$
Thus, by Proposition 5.2 $ii)$ and the fact $\b_{\lambda'}-\f_{\lambda'} \in v\mathcal{L}$, we have
$$\b_\m|0,\s_l\rangle \equiv(-v)^{(n-1)|\mu|}a_\mu \f_{\lambda'}|0,\s_l\rangle \equiv (-v)^{(n-1)|\mu|}a_\mu
\b_{\lambda'}|0,\s_l\rangle\;(\mathrm{mod}\;\mathcal{L}_{\Lambda_{\s_l}})$$
Now by Proposition 5.2 $i)$,
$$(-v)^{(n-1)|\mu|}a_\mu\b_{\lambda'}|0,\s_l\rangle \equiv |\nu_l,\s_l\rangle$$
for some $\nu_l \in \Pi^l$ if $\b_{\lambda'}|0,\s_l\rangle \in \mathbf{B}^+_{\s_l}$ and
$(-v)^{(n-1)|\mu|}a_\mu\b_{\lambda'}|0,\s_l\rangle =0$
if $\b_{\lambda'}|0,\s_l\rangle=0$. Moreover, by \cite{U}, Proposition 4.12, $\overline{\b_\m|0,\s_l\rangle}=\overline{\b_\m}
\;\overline{|0,\s_l\rangle}=\b_\m|0,\s_l\rangle$. Hence $\b_\m|0,\s_l\rangle=\b^+_{\nu_l,\s_l}$ in the first case and
$\b_\m|0,\s_l\rangle =0$ in the second case. \qed
\paragraph{}Let $\mathcal{C}^0_{\s_l}$ be the subgraph of the crystal graph of $\Lambda_{\s_l}$ corresponding to $\U^-_n |0,\s_l\rangle$. The above theorem implies the following special case of the positivity conjecture of Uglov (see \cite{U}, Section 4):
\begin{cor} For any $\s_l \in \Z^l(s)$ and any $\lambda_l \in \mathcal{C}^0_{\s_l}$ we have
$$\b_{\lambda_l,\s_l} \in \bigoplus_{\mu_l} \N[v] |\mu_l,\s_l\rangle.$$
\end{cor}
\noindent
\section{Appendix}
\paragraph{6.1}\textit{Proof of Lemma 2.2.} Let us prove that for any $u \in \uan$ we have :
$$[\rho(e_i),\rho(u)]=\rho\big(\frac{\k_ie_i''(u)-\k_i^{-1}e'_i(u)}{v-v^{-1}}\big).$$
This result is well-known for $\mathbf{U}^-(\hat{\mathfrak{sl}}_n)$ (see \cite{K}, Lemma 3.4.1). Let us prove it for $\f_{k\s}, l \in \N^*$. We use the presentation of $\La$ in terms of Young diagrams and the description of the representation $\rho$ in terms of the Hall algebra $\mathbf{U}^-_{\infty}$ of the infinite quiver and the quantum enveloping algebra $\mathbf{U}(\mathfrak{sl}_\infty)$, as given in \cite{VV}. We keep the notations of \cite{VV}, Section 6. In particular, let $\gamma_d: \uan \to \mathbf{U}^-_{\infty}$ be the map defined in \cite{VV}, Section 6; For any $\mathbf{d}=(d_1,d_2,\ldots) \in \N^{(\Z)}$ and $k \in \N$ we put
$$
h(\mathbf{d})=\sum_{i<j,i \equiv j} d_i(d_{j+1}-d_j), \qquad \mathbf{d}'=\sum_{j<m;j \equiv\;m} d_j \epsilon_m, \qquad \mathbf{k}''=\sum_{r<k;r \equiv\;k} \epsilon_r.$$
Let $\pi: \N^{(\Z)} \to \N^{I}$ be the reduction modulo $n$. To avoid confusion, we will denote by $\underline{\e}_i,\;\underline{\f}_i$ the elements of $\mathbf{U}(\hat{\mathfrak{sl}}_n)$ and $\uan$, and by $\e_i,\f_i$ the elements of $\mathbf{U}(\hat{\mathfrak{sl}}_\infty)$. Finally, to simplify notations, we will write $\sum_{\alpha<\beta}$ for $\sum_{\alpha<\beta;\alpha \equiv \beta}$. We have, in $\mathrm{End}(\Lambda^\infty)$,
\begin{align*}
\underline{\f}_{k\s}&=\sum_{\d \in \pi^{-1}(k\s)} v^{h(\d)}\f_\d \k_{\d'},\\
\underline{\e}_i&=\sum_{j \equiv i} \e_j\k_{\mathbf{j}''}^{-1}.
\end{align*}
 By \cite{VV}, Section 5.2, the element $\f_\d |\lambda\rangle$ is zero for all $\lambda$ if there exists $j$ such that $d_j \not\in \{0,1\}$. Hence,
\begin{align*}
[\underline{\e}_i,\underline{\f}_{k\s}]&=\sum_{\underset{\d}{j \equiv i}} v^{h(\d)}[\e_j \k_{j''}^{-1},\f_\d \k_{\d'}]\\
&=\sum_{\underset{\d}{j \equiv i}} v^{h(\d)}\big(v^{\sum_{l<j} \text{wt}_l(\d)} \e_j\f_\d-v^{-\sum_{t<m} \d_t \text{wt}_m(-\epsilon_j)}\f_\d \e_j\big) \k_{j''}^{-1}\k_{\d'}\\
&=\sum_{\underset{\d}{j \equiv i}} v^{h(\d)+\sum_{l<j} 2d_l-d_{l-1}-d_{l+1}} [\e_j,\f_\d]  \k_{j''}^{-1}\k_{\d'}
\end{align*}
Recall that $\f_\d$ is equal to the ordered product $\f_\d=\ldots \f_{i-1}^{d_{i-1}} \f_{i}^{d_{i}} \f_{i+1}^{d_{i+1}}\ldots$. Then
$$[\e_j,\f_\d]=\cases \left(\prod_{r=-\infty}^{j-1} \f_r^{d_r}\right)\frac{\k_j-\k_j^{-1}}{v-v^{-1}}  \left(\prod_{s=j+1}^{\infty} \f_s^{d_s}\right) &if\; d_j=1\\
&\\
0&if\; d_j=0 \endcases.$$
Thus,
\begin{align*}
[\underline{\e}_i,\underline{\f}_{k\s}]&= \sum_{\underset{\d\;s.t\; d_j=1}{j\equiv i}} v^{h(\d)+\sum_{l<j} 2d_l-d_{l-1}-d_{l+1}}\left(\prod_{r=-\infty}^{j-1} \f_r^{d_r}\right)\frac{\k_j-\k_j^{-1}}{v-v^{-1}} \left(\prod_{s=j+1}^{\infty} \f_s^{d_s}\right) \k_{j''}^{-1}\k_{\d'}\\
&=\sum_{\underset{\d\;s.t\; d_j=1}{j\equiv i}} v^{h(\d)+\sum_{l<j} 2d_l-d_{l-1}-d_{l+1}}\left(\prod_{r=-\infty}^{j-1} \f_r^{d_r}\right) \left(\prod_{s=j+1}^{\infty} \f_s^{d_s}\right)\frac{v^{d_{j+1}}\k_j-v^{-d_{j+1}}\k_j^{-1}}{v-v^{-1}} \k_{j''}^{-1}\k_{\d'}.
\end{align*}
We have
$$h(\d) +\sum_{l<j} (2d_l-d_{l-1}-d_{l+1})=h(\d-\epsilon_j)+ \sum_{s>j} (d_{s+1}-d_s)+\sum_{l<j} (d_l -d_{l+1})$$
and
 $$\k_{j''}^{-1}\k_{\d'}=\k_{(\d-\epsilon_j)'}\prod_{t<j}\k_t^{-1}\prod_{r>j}\k_r.$$
Therefore
$$
[\underline{\e}_i,\underline{\f}_{k\s}]=\sum_{\tilde{\d} \in \pi^{-1}(k\s-\epsilon_i)} v^{h(\tilde{\d})}\f_{\tilde{\d}} \k_{\tilde{\d}'} \sum_{\underset{\tilde{\d}_j=0}{j \equiv i}} v^{\sum_{s>j}(d_{s+1}-d_s)+\sum_{l<j}(d_l-d_{l+1})} \left(\frac{v^{d_{j+1}}\k_j-v^{-d_{j+1}}\k_j^{-1}}{v-v^{-1}} \right)\prod_{t<j}\k_t^{-1}\prod_{r>j}\k_r.$$
Moreover, $\f_{\tilde{\d}}|\lambda\rangle\neq 0$ if and only if $\k_t |\lambda\rangle=v |\lambda\rangle$ for all $t$ such that $\tilde{d}_t=1$. It follows that, in $\mathrm{End}(\Lambda^\infty)$,
\begin{align*}
[\underline{\e}_i,\underline{\f}_{k\s}]&=\sum_{\tilde{\d} \in \pi^{-1}(k\s-\epsilon_i)} v^{h(\tilde{\d})}\f_{\tilde{\d}} \k_{\tilde{\d}'} \sum_{\underset{\tilde{\d}_j=0}{j \equiv i}} v^{\sum_{s>j}d_{s+1}-\sum_{l<j}d_{l+1}} \left(\frac{v^{d_{j+1}}\k_j-v^{-d_{j+1}}\k_j^{-1}}{v-v^{-1}} \right)\prod_{\underset{\tilde{d}_t=0}{t<j}}\k_t^{-1}\prod_{\underset{\tilde{d}_r=0}{r>j}}\k_r\\
&=\sum_{\tilde{\d} \in \pi^{-1}(k\s-\epsilon_i)} v^{h(\tilde{\d})}\f_{\tilde{\d}} \k_{\tilde{\d}'}\left( \prod_{t \equiv i} \k_t-\prod_{r \equiv i} \k_r^{-1}\right)
\end{align*}
Note that both infinite products make sense as elements of $\mathrm{End}(\Lambda^\infty)$ since for any $\lambda \in \Pi$, $\k_t|\lambda\rangle=|\lambda \rangle$ for all but a finite number of values of $t$. Hence
\begin{align*}
[\underline{\e}_i,\underline{\f}_{k\s}]
&=\underline{\f}_{k\s-\epsilon_i}\frac{v\k_i-v^{-1}\k_i^{-1}}{v-v^{-1}}\\
&=\frac{v^{-1}\k_i-v\k_i^{-1}}{v-v^{-1}}\underline{\f}_{k\s-\epsilon_i}
\end{align*}
as desired.\\
Now suppose that Lemma 2.2 holds for $u \in \uan[\alpha]$ and $w \in \uan[\beta]$. Then
\begin{equation*}
\begin{split} 
\rho&(\underline{\e}_i)\rho(uw)\\
&=\rho(u)\rho(\underline{\e}_i)\rho(w)+\rho\left(\frac{\k_ie''_i(u)-\k_i^{-1}e'_i(u)}{v-v^{-1}}w\right)\\
&=\rho(uw)\rho(\underline{\e}_i)+\rho\left(u\frac{\k_ie_i''(w)-\k_i^{-1}e_i'(w)}{v-v^{-1}}\right) +\rho\left(\frac{\k_ie''_i(u)-\k_i^{-1}e'_i(u)}{v-v^{-1}}w\right).
\end{split}
\end{equation*}
Therefore
\begin{equation*}
\begin{split}
[\rho&(\underline{\e}_i),\rho(uw)]\\
&=\rho\left(\frac{\k_i(v^{\text{wt}_i(\alpha)}ue''_i(w)+e''_i(u)w)-\k_i^{-1}(v^{-\text{wt}_i(\alpha)}ue'_i(w)+e'_i(u)w)}{v-v^{-1}}\right)\\
&=\rho\left(\frac{\k_ie''(uw)-\k_i^{-1}e'_i(uw)}{v-v^{-1}}\right),
\end{split}
\end{equation*}
i.e. Lemma 2.2 is true for the product $uw$. This proves Lemma 2.2 since, by \cite{GP} Theorem 3.1, $\ukn$ is generated by $\mathbf{U}^-(\hat{\mathfrak{sl}}_n)$ and the elements $\f_{l\s},\; l \in \N^*$.\qed
\paragraph{6.2}\textit{Proof of Lemma 3.4.} Fix integers $k,l$ with $k>l$, and let $\mu \in \Pi$, $|\mu|=k-l$. Fix a $I$-graded vector space $V$ of dimension $k\s$, a subspace $V' \subset V$ of dimension $(k-l)\s$ and choose an element $y \in O_{\mathbf{m}(n\mu)} \subset E_{V'}$. Finally, let $\m \in \mathcal{M}$ with $\mathrm{dim}\;\m=k\s$. By definition, 
$$\Delta_{(1)^l,\mu}(\mathbf{1}_\m)(\f_{l\s} \otimes \f_{n\mu})=\# I^\m_{|q=v^{-2}} \mathbf{1}_{l\s} \otimes \mathbf{1}_{n\mu},$$
where $I^\m=\{x \in O_\m \subset E_V\;|\;x_{|V'}=y,\;\mathrm{Im}\;x \subset V'\}$. Recall that $\#I^\m$ is a polynomial in $q$, and that by the Lang-Weil theorem (\cite{LW}) we have $\#I^\m=q^{\mathrm{dim}\;I^\m}+ \mathcal{O}(q^{\mathrm{dim}\;I^\m-1})$. Moreover it follows from \cite{McD}, (6.17), that $c^\lambda_{(1)^l,\mu} \in \{0,1\}$ and that $I^\m=\emptyset$ if $\m=\m(n\lambda)$ with $c^{\lambda}_{(1)^l,\mu}=0$. Hence Lemma 3.4 is equivalent to the following dimension inequalities :
\begin{equation}\label{E:60}
\begin{split}
\mathrm{dim}\;O_\m-\mathrm{dim}\;O_{n\mu} &\geq 2\;\mathrm{dim}\; I^\m,\\
\mathrm{dim}\;O_\m-\mathrm{dim}\;O_{n\mu}&= 2\;\mathrm{dim} \;I^\m \Leftrightarrow\; \m=\m(n\lambda)\;\mathrm{for\;some\;}\lambda \in \Pi\;\mathrm{with}\; c^\lambda_{(1)^l,\mu}=1.
\end{split}
\end{equation}
For $i \in \N^*$, $j \in I$ and for any $x \in O_\m$, let us set
\begin{align*}
d^i_j=&\mathrm{dim\;Ker\;}y^i_{|V'_j}-\mathrm{dim\;Ker\;}y^{i-1}_{|V'_j},\\
\td^i_j=&\mathrm{dim\;Ker\;}x^i_{|V_j}-\mathrm{dim\;Ker\;}x^{i-1}_{|V_j},\\
\theta^i_j=&\td^i_j-d^i_j.
\end{align*}
Then $I^\m \neq \emptyset$ if and only if $\theta^i_j \in \N$ and $\sum_i \theta^i_j=l$ for all $j \in I$. Note that $d^k_j=d^k_{j'}$ for any $k,j,j'$ since $\m(n\mu) \in \mathcal{M}^{\mathrm{per}}.$ A direct computation gives
\begin{equation*}
\mathrm{dim}\; O_\m=\sum_j \left( \sum_{k>i} (\td^k_j+\td^k_{j-1})\td^i_j\right)=\sum_j \left( \sum_{k>i} (d^k_j+d^k_{j-1}+\theta^k_j+\theta^k_{j-1})(d^i_j+\theta^i_j)\right).
\end{equation*}
Similarly, 
$$\mathrm{dim}\; O_{\m(n\mu)}=\sum_j \left( \sum_{k>i} (d^k_j+d^k_{j-1})d^i_j\right).$$
Thus,
\begin{equation}\label{E:62}
\mathrm{dim}\;O_\m-\mathrm{dim}\;O_{\m(n\mu)}=\sum_j \left(2\sum_{k \neq i} \theta^k_j d^i_j + \sum_{k>i}(\theta^k_j + \theta^k_{j-1})\theta^i_j\right).
\end{equation}
Now we compute $\mathrm{dim}\;I^\m$. Fix a complementary subspace $U$ of $V'$ in $V$. An element $x \in E_V$ satisfying $\mathrm{Im}\;x \subset V'$ and $x_{|V'}=y$ is uniquely determined by the collection of maps $x: U_j \to V'_{j+1}$. Moreover, 
\begin{equation*}
\begin{split}
\mathrm{dim}\;(\mathrm{Ker}\;x^k_{|V_j})&=\mathrm{dim}\;(\mathrm{Ker}\;
y^k_{|V'_j}) +\mathrm{dim}\;\mathrm{Ker}\;x^k_{|U_j}+\mathrm{dim}\; (x^k(U_j) \cap y^k(V'_j))\\
&=\mathrm{dim}\;(\mathrm{Ker}\;y^k_{|V'_j}) +\mathrm{dim\;Ker}\;x_{|U_j}+\mathrm{dim}\;(x(U_j) \cap(y^{k-1})^{-1}(y^k(V'_j))).
\end{split}
\end{equation*}
Set $Y^1_k=\{0\}$ and $Y^k_j=(y^{k-1})^{-1}(y^k(V'_j))$ if $k>1$. Then, for $k>1$ 
\begin{equation}\label{E:621}
\mathrm{dim}\; Y^k_{j}=\sum_{i<k} d^i_{j+1} + \sum_{i>k} d^i_j=\sum_{i \neq k} d^i_j.
\end{equation}
Now,
\begin{equation*}
\mathrm{dim}\;(\mathrm{Ker}\;x^k_{|V_j})-\mathrm{dim}\;(\mathrm{Ker}\;x^{k-1}_{|V_j})
=d^k_j+ \mathrm{dim}(x(U_j)\cap Y^k_j) - \mathrm{dim}(x(U_j)\cap Y^{k-1}_j).
\end{equation*}
Hence $x \in O_\m$ if and only if for all $j \in I$
\begin{align}
\mathrm{dim\;Ker}\;x_{|U_j}&=\theta^1_k,\notag\\
\forall \;k>1,\qquad  \mathrm{dim}(x(U_j)\cap Y^k_j) - &\mathrm{dim}(x(U_j)\cap Y^{k-1}_j)=\theta^k_j.\label{E:622}
\end{align}
The variety $X_j$ of subspaces $W_j \subset V'_{j+1}$ satisfying (\ref{E:622}) is of dimension
$$\mathrm{dim}\;X_j=\sum_{l\geq 2} (\mathrm{dim}\;Y^l_j-\sum_{t=2}^l\theta^t_j)\theta^l_j.$$
A direct computation using (\ref{E:621}) now gives
\begin{equation}\label{E:623}
\mathrm{dim}\;I^\m=\sum_j \left(\sum_{k \neq i} \theta^k_jd^i_j + \sum_{k>j} \theta^k_j \theta^i_j\right)
\end{equation}
Thus the dimension inequalities (\ref{E:60}) are consequences of
(\ref{E:62}), (\ref{E:623}) and the following result :\\ 
\textit{Claim.} For all collections of positive integers $(\theta^i_j)_{i,j}$, $i=1,\ldots h$, $j \in \Z/n\Z$ satisfying $\sum_i \theta^i_j=l$ there holds
\begin{equation}\label{E:claim}
\sum_j\left(\sum_{k>i} \theta^k_{j-1}\theta^i_j\right) \geq \sum_j \left(\sum_{k>i} \theta^k_j \theta^i_j\right),
\end{equation}
with equality if and only if $\theta^k_j=\theta^k_{j'}$ for all $k,j,j'$.\\
\hbox to1em{\hfill}We argue by induction on $h$ to prove (\ref{E:claim}). The claim is trivial if $h=1$. Suppose that (\ref{E:claim}) is proved for all $h'<h$. We first note that
\begin{equation}\label{E:cl1}
\begin{split}
\sum_j \sum_{k>i} \left( \theta^k_{j-1}\theta^i_j-\theta^k_j\theta^i_j\right)=& \sum_j \left(\sum_{k>i}^{h-1} \theta^k_{j-1}\theta^i_j-\theta^k_j\theta^i_j + \sum_{i<h} (l-\sum_{k<h}\theta^k_j)(\theta^i_{j-1}-\theta^i_j)\right)\\
=&\sum_j\sum_{k\geq i}^{h-1} \theta^k_j(\theta^i_j-\theta^i_{j-1}).
\end{split}
\end{equation}
Let us freeze variables $\theta^k_j$, $k>1$ and consider $G((\theta^1_j))=\sum_j\sum_{k\geq i}^{h-1} \theta^k_j(\theta^i_j-\theta^i_{j-1})$. Then a direct computation shows that $G((\theta^1_j))$ reaches its global minimum when for all $j$
$$(\theta^1_j-\theta^1_{j-1})-(\theta^1_{j+1}-\theta^1_j)=\sum_{k=2}^{h-1} (\theta^k_{j+1} + \theta^k_{j-1} -2 \theta^k_j),$$
i.e. when for all $j$
$$\sum_{k=1}^{h-1}(\theta^k_{j+1}- \theta^k_j)=\sum_{k=1}^{h-1}(\theta^k_{j}- \theta^k_{j-1}).$$
Since $\sum_j\sum_{k=1}^{h-1}(\theta^k_{j+1}- \theta^k_j)=0$, this implies that $\sum_{k=1}^{h-1}\theta^k_{j+1}=\sum_{k=1}^{h-1} \theta^k_j$ for all $j$. But then $\theta^h_j=l-\sum_{k=1}^{h-1}\theta^k_j=\theta^h_{j+1}$, and in this case
$$ 
\sum_j\sum_{k>i}^h (\theta^k_{j-1}-\theta^k_j)\theta^i_j=\sum_j\sum_{k>i}^{h-1} (\theta^k_{j-1}-\theta^k_j)\theta^i_j.$$
The result then follows from the induction hypothesis.\qed

\vspace{.2in}
{\centerline{\textbf{Acknowledgments}}}
\paragraph{}I would like to thank my advisor Eric Vasserot for
suggesting this problem to me and for his patience and guidance. I am
endebted to B. Leclerc for many enlightning discussions and comments
on this paper and to A. Braverman and D. Gaitsgory for interesting discussions.

\small{
\vspace{4mm}
Olivier Schiffmann, ENS Paris, 45 rue d'Ulm, 75005
PARIS; \texttt{schiffma@clipper.ens.fr}
\end{document}